  \documentstyle[12pt]{article}
\newcommand{\mylabel}[1]{
\label{#1}}

\newcommand{\numero}[1]{
\addtocounter{section}{1}
\begin{center}{\bf \thesection .\
#1\vspace{-.1in}}\end{center}
\setcounter{subsection}{0}
\setcounter{lemma}{0}\indent}

\newcommand{\subnumero}[1]{
\pagebreak[1]\begin{center}{\em #1}\nopagebreak\end{center}
}

\newcommand{\eop}{\hfill $\Box$\vspace{.1in}}

\newtheorem{lemma}{Lemma}[section]
\newtheorem{theorem}[lemma]{Theorem}

\newtheorem{corollary}[lemma]{Corollary}

\newtheorem{proposition}[lemma]{Proposition}

\newcommand{\cc}{{\bf C}}

\newcommand{\qq}{{\bf Q}}

\newcommand{\zz}{{\bf Z}}
\newcommand{\pp}{{\bf P}}

\newcommand{\Cc}{{\cal C}}

\newcommand{\Ff}{{\cal F}}
\newcommand{\Gg}{{\cal G}}

\newcommand{\Oo}{{\cal O}}

\newcommand{\Hh}{{\cal H}}

\newcommand{\Xx}{{\cal X}}

\newcommand{\Gm}{{\bf G}_m}

\newcommand{\gerg}{{\bf g}}

\newcommand{\delbar}{\overline{\partial}}

\begin{document}

\section*{Secondary Kodaira-Spencer classes and nonabelian Dolbeault cohomology}

Carlos Simpson\newline
CNRS, Laboratoire Emile Picard\newline
Universit\'e Toulouse III\newline
31062 Toulouse CEDEX, France\newline
carlos@picard.ups-tlse.fr

\bigskip

One of the nicest things about variations of Hodge structure is the
``infinitesimal variation of Hodge structure'' point of view \cite{IVHS}.
A variation of Hodge structure $(V=\bigoplus V^{p,q}, \nabla )$
over a base $S$
gives rise at any point $s\in S$ to the {\em Kodaira-Spencer map}
$$
\kappa _s: T(S)_s \rightarrow Hom (V^{p,q}_s, V^{p-1, q+1}_s).
$$
In the geometric situation of a family $X\rightarrow S$ we have
$$
V^{p,q}_s = H^q(X_s, \Omega ^p_{X_s}),
$$
and the Kodaira-Spencer map is given by cup-product with the Kodaira-Spencer
deformation class
$$
T(S)_s \rightarrow H^1(X_s, T(X_s)).
$$
The Kodaira-Spencer map is a component of the connection $\nabla$. In
particular, this implies that if $\kappa _s\neq 0$ then the connection
$\nabla$ is nontrivial with respect to the Hodge decomposition. Various
Hodge-theory facts imply that the global monodromy must be nontrivial in this
case. We can be a bit more precise: if $u\in V^{p,q}$ is a vector such that
$\kappa _s(v)(u)\neq 0$ for some tangent vector $v\in T(S)_s$, then $u$ cannot
be preserved by the global monodromy. Thus a local calculation (which actually
only depends on the first-order deformation of $X_s$) implies a global fact. In
particular this global fact would hold for any family of varieties $X'$ over
any base $S'$, such that the new family osculates to order $1$ with the
original one (say as a map from $S'$ into the moduli stack of the fibers).
A particularly nice aspect of this situation is that the Kodaira-Spencer map is
defined on the Dolbeault cohomology $H^q(X_s, \Omega ^p)$ and in particular it
is obtained involving only algebraic-geometric calculations (just a cup-product
with the deformation class)---no analytic considerations are needed.

The goal of this paper is to calculate an example showing a similar type of
behavior with a secondary Kodaira-Spencer class coming from nonabelian
cohomology with coefficients in the complexified $2$-sphere $T=S^2 \otimes \cc$.
For such a $T$ (or any other coefficient stack similar in nature) we define
the {\em nonabelian Dolbeault cohomology of $X$ with coefficients in $T$},
denoted $Hom (X_{Dol}, T)$.  When $X$ varies in a family parametrized by a base
$S$ then we show how to define a secondary class which is a map
$$
\bigwedge ^2T(S)_s\rightarrow
\pi _1(Hom (X_{Dol}, T)) = H^2_{Dol}(X)/(\eta )
$$
where $\eta$ is the class in $H^2_{Dol}(X)$ pulled back from the tautological
class on $T$ by the map $X_{Dol}\rightarrow T$ which we take as basepoint.
The secondary class is defined when the primary Kodaira-Spencer classes vanish.

It seems likely, although we don't show that here, that
our secondary class is just a quadruple Massey
product
$$
\alpha \wedge \beta \mapsto \{ \eta , \eta , \alpha , \beta \}
$$
for $\alpha , \beta \in H^1(X, TX)$ and $\eta \in H^2_{Dol}(X)$ with $\eta \cup
\eta = 0$. Instead of interpreting the secondary class this way, we calculate
it directly using Mayer-Vietoris arguments for nonabelian cohomology and
reducing to calculations in abelian cohomology.

Our
main purpose in the present paper is to calculate the secondary class in a
specific, somewhat instructive, example. We look at the family of varieties $X$
obtained by blowing up a point $P\in Z$ on a smooth surface $Z$. This family is
parametrized in an obvious way by $Z$. If $Z$ is simply connected, then any
standard Hodge-theoretic information related to $X$ must be independent of the
basepoint because it would vary in a variation of (mixed) Hodge structure
parametrized by $Z$, and such a variation is forcibly constant. However, we will
show that the secondary Kodaira-Spencer class for nonabelian Dolbeault
cohomology is nonzero, if $Z$ has a nonzero holomorphic $2$-form.

Before getting to the example, we will discuss some general aspects of
nonabelian cohomology with emphasis on the case of Dolbeault cohomology. We
start with a quick review of $n$-stacks in \S 1. Then we define the notion of
``connected very presentable shape''in \S 2 by looking at maps to $n$-stacks $T$
which are $0$-connected, with $\pi _1(T)$ an affine algebraic group scheme, and
$\pi _i(T)$ vector spaces, for $i\geq 2$. The next step in \S 3 is to set
$X_{Dol}$ equal to the $1$-stack over $X$ whose fiber is $K(\widehat{TX}, 1)$
where $\widehat{TX}$ is the completion of the tangent bundle of $X$ along the
zero-section.  We define the {\em nonabelian Dolbeault cohomology of $X$ with
coefficients in an $n$-stack $T$} to be $Hom (X_{Dol},T)$. We look at this
particularly for connected very presentable $T$. This gives rise to the {\em
Dolbeault shape} which is the $(n+1)$-functor
$$
T\mapsto Hom (X_{Dol}, T)
$$
on connected very presentable $T$. Sections  4 and 5 treat group actions and
secondary classes in the $n$-stack situation. In \S 5 we define the secondary
Kodaira-Spencer class in general. Then in \S 6 we look at a particular case of
$T$, namely the {\em complexified $2$-sphere}, defined by the conditions that
$\pi _2=\pi _3 = \Oo $ and having nontrivial Whitehead product.

With these preliminary steps done, we get to our example in \S\S 7-8. The main
result is Theorem \ref{calculation}. It says that if $X$ is the blowup of a
surface $Z$ at a point $P$ then the secondary Kodaira-Spencer class for
nonabelian Dolbeault cohomology $Hom (X_{Dol}, T)$ with coefficients in
the complexified $2$-sphere $T$, is Serre dual to the evaluation map of
holomorphic $2$-forms at $P$. In particular if $H^0(Z,\Omega ^2_Z)\neq 0$ then
this class is nonzero.

At the end in two appendices, we will discuss some topics from \cite{kobe} but
in greater detail. First, the ``Breen calculations'' giving the cohomology of
$K(\Oo , n)$ or more generally $K(V/S,n)$ for a vector sheaf $V$ over a base
scheme $S$. Then we discuss representability of certain shapes.
These topics come up in a few
places in the body of the paper, which is the reason for the appendices; for
an introduction we refer the
reader to the appendices.

Without going into great bibliographic detail, we point out here some recent
papers which seem to be somewhat related. Karpishpan \cite{Karpishpan1}
\cite{Karpishpan2} treats higher-order Kodaira-Spencer mappings, i.e. the
higher order derivatives of the period map. The same is treated by Ran
\cite{Ran}
and Esnault-Viehweg \cite{EsnaultViehweg}.  Biswas
defines secondary invariants for families of Higgs bundles \cite{Biswas}. Also
Bloch and Esnault treat algebraic Cherns-Simons classes \cite{BlochEsnault},
which are types of secondary classes.

I would like to thank Mark Green for an inspiring question---albeit one which
the present paper doesn't answer. He asked whether there are examples of
families of varieties where the variation of Hodge structure on the
cohomology is constant, but where the variation of mixed Hodge structure on the
homotopy groups is nontrivial. In the absence of an answer to that question
(which is very interesting), the typical mathematician's reply is to change the
question---in this case, to look for an example where even the mixed Hodge
structures on the homotopy groups remain constant, but where the Hodge
filtration on the full homotopy type is nonconstant.
\footnote{
After the first version of this paper, Richard Hain pointed out to me the
following example answering Mark Green's original question. The example comes
from  a paper of Carlson, Clemens, and Morgan (\cite{CarlsonEtAl}, p. 330). Let
$C\subset \pp ^3$ be an embedded curve of positive genus. For points $p,q\in C$,
let $X_{p,q}$ be the $3$-fold obtained by first blowing up $p$ and $q$ and then
blowing up the strict transform of $C$.  The family of $X_{p,q}$ parametrized by
(an open subset of) $C\times C$ has constant variation of Hodge structures on
the  cohomology but, according to \cite{CarlsonEtAl} the variation of MHS on the
homotopy is nonconstant.}

\subnumero{Notation}

We always work in characteristic $0$. In order to simplify notation we use
$\cc$ as the ground field (i.e. $Spec (\cc )$ as base scheme), but everything we
say would work equally well over any ground field of characteristic $0$. Let
$Sch /\cc$ denote the site of schemes of finite type over $Spec (\cc )$ with
the etale topology.

The structure sheaf $\Oo$ on $Sch /\cc$ is the sheaf defined by
$$
\Oo (Y):= \Gamma (Y, \Oo _Y).
$$
It is represented by the affine line ${\bf A}^1$, in other words it is
represented by the $1$-dimensional vector space $\cc$. A finite dimensional
vector space represents a sheaf of the form $\Oo ^a$.

\numero{Basic remarks about $n$-stacks}

We make some brief remarks about $n$-stacks as we shall use them in this paper.
For all details the reader is referred to the following references:
\newline
---for the history and basic notions of simplicial presheaves: Brown
\cite{Brown}, Illusie \cite{Illusie}, Jardine \cite{Jardine};
\newline
---for the history and basic notions of $1$-stacks: Artin
\cite{ArtinInventiones}, Deligne-Mumford \cite{Deligne-Mumford}, and
Laumon-Moret-Bailly \cite{LMB};
\newline
---for cohomological theory using simplicial presheaves: Thomason
\cite{Thomason};
\newline
---for a ``homotopy coherent'' approach: Cordier-Porter \cite{CordierPorter}
and also \cite{flexible};
\newline
---for $n$-categories and $n$-stacks: Grothendieck \cite{PursuingStacks},
Breen \cite{BreenAsterisque}, Gordon-Power-Street
\cite{GordonPowerStreet}, Tamsamani \cite{Tamsamani}, Baez-Dolan
\cite{Baez-Dolan}, and several papers of the author.

The first main remark is that we shall almost always be concerned with
$n$-stacks of $n$-groupoids, and following the intuition put forth in
\cite{PursuingStacks}, an $n$-groupoid is the same thing (up to homotopy) as a
topological space whose homotopy groups vanish in degrees $i>n$---we call such
a space {\em $n$-truncated}. Thus it is safe to replace $n$-groupoids
everywhere by $n$-truncated topological spaces or, again equivalently,
$n$-truncated simplicial sets. In this point of view, an $n$-stack (on the
site $Sch /\cc $ which is fixed throughout) is just a presheaf of
simplicial sets
otherwise known as a {\em simplicial presheaf}, which is object-by-object
$n$-truncated. Once one has made the passage to simplicial presheaves, the
truncation condition is no longer crucial (although it often facilitates
arguments and many things in the literature are only stated in this case or
under a complementary hypothesis about vanishing cohomological dimension).
Thus, when we speak of ``$n$-stacks'', one way to read this is in terms of the
homotopy theory of simplicial presheaves, see \cite{Brown} \cite{Illusie}
\cite{Jardine} \cite{Thomason}.

Another reading would plunge directly into
the theory of $n$-categories and $n$-stacks, imposing the groupoid condition
along the way. For $n\leq 3$ (which in the end is the case we treat in the
present paper) this can be had in a relatively formulaic way in
\cite{BreenAsterisque} and \cite{GordonPowerStreet}. For arbitrary $n$, see
\cite{Tamsamani} \cite{nCAT}, but unfortunately the $n$-stack part of this
theory still needs to be worked out a bit more.

The only place where we make reference to $n$-categories which are not
$n$-groupoids is when we look at the $n+1$-category $nSTACK$ of $n$-stacks (of
groupoids). This $n+1$-category has the property of being {\em $1$-groupic},
i.e. the morphism $n$-categories are $n$-groupoids. One can safely replace the
morphism $n$-groupoids by spaces or simplicial sets, and one obtains the
notion of {\em Segal category} \cite{effective}, motivated by Segal's delooping
machine \cite{Segal}. This notion came into higher category theory in
Tamsamani's definition of $n$-category \cite{Tamsamani}.  Thus $nSTACK$ may be
considered as a Segal category. This fits in relatively nicely with the
simplicial presheaf point of view; in fact this Segal category comes from the
simplicial category of fibrant and cofibrant objects in Jardine's closed model
category of simplicial presheaves.  This comes up in
looking at the functoriality in $T$ of the construction $Hom (X_{Dol}, T)$.
In one other place we refer to the $n+1$-stack $n\underline{STACK}$ of
$n$-stacks, which is discussed somewhat in
\cite{nCAT} and \cite{limits}; we don't get any further into the general theory
here.

Currently, the
``simplicial presheaves'' alternative is the most accessible (it is also
historically the first, dating from \cite{Brown}). In terms of simplicial
presheaves, an {\em $n$-stack} is a simplicial presheaf on $Sch/\cc $ which is
object-by-object $n$-truncated. If $X,Y$ are simplicial presheaves, then we
obtain a simplicial presheaf $Hom(X,Y)$ by first replacing $Y$ by a fibrant
object \cite{Jardine}, then looking at the internal $Hom$ of simplicial
presheaves. In particular this gives a simplicial set $Hom(X,Y)(Spec (\cc ))$
and this family of simplicial sets makes the simplicial presheaves into a
simplicial category (or Segal category) which we denote $nSTACK$. When we speak
of morphisms between $n$-stacks, the above procedure is always understood, i.e.
we always replace the target by a fibrant object.

Given a presheaf of $n$-groupoids or presheaf of spaces, we often want to take
the ``associated $n$-stack''. What this means depends somewhat on the point of
view which is taken. If one works with objects in a closed model category such
as that of simplicial presheaves, then this just means to consider the object
as an element of the closed model category. One might also want to say that it
means to replace the object by a weakly equivalent fibrant object. Finally
there is an intermediate notion based on enforcing the global descent condition
but not the local fibrant condition. It doesn't really matter which point of
view we adopt, since when looking at morphisms to a given object, we always
replace it by an equivalent fibrant object anyway.

The ``yoga'' of the situation is that one can do topology with $n$-stacks
instead of spaces (or more precisely $n$-truncated spaces).  In particular, all
standard constructions and results in  algebraic topology carry over to
$n$-stacks. Most of these are contained somewhere in the literature referred to
above; but if not, we don't give proofs here as that would get beyond the scope
of the present paper.

There is basically only one slight ``twist'' which is not present in the
topological case: this is that the $0$-truncated objects  can be topologically
nontrivial, i.e. can have cohomology. In the usual topological case, the
$0$-truncated objects are just the disjoint unions of contractible components
and these make no significant contribution to homotopy. In the case of
$n$-stacks over a site, one can have $0$-stacks, i.e. sheaves of sets, with
nontrivial cohomology (this is the case of a smooth projective variety $X$, for
example); and similarly there are sheaves of groups over $\ast$ (the site
itself) which can have nontrivial cohomology. The upshot of all this is that
when it comes to choosing basepoints for an $n$-stack $T$, one must choose first
an object $Y\in Sch /\cc$ and then choose a basepoint $t\in T(Y)$. (In other
words, if we look only at basepoints in $T(Spec \cc )$ we might be missing
some topology.)

The first place where the previous paragraph has an impact is in the notion of
{\em homotopy groups}. If $T$ is an $n$-stack then for any
$Y\in Sch /\cc$ and $t\in T(Y)$ we obtain a presheaf, denoted in utmost
precision by
$$
\pi _i^{\rm pre}(T|_{Sch /Y}, t)
$$
but which we often shorten to $\pi _i^{\rm pre}(T,t)$. This is a presheaf of
groups over the site $Sch /Y$, abelian if $i\geq 2$. On the other hand the
presheaf $\pi _0^{\rm pre}(T)$ is defined absolutely as a presheaf of sets over
$Sch /\cc $.

The definition which is fundamental to the theory is that we define
$$
\pi _i(T|_{Sch /Y}, t)
$$
to be the sheaf associated to the presheaf
$\pi _i^{\rm pre}(T|_{Sch /Y}, t)$. Similarly $\pi _0(T)$ is the sheaf of sets
associated to the presheaf $\pi _0^{\rm pre}(T)$. These, and not the presheaf
versions, are the only thing we care about. This is formalized by saying that
a morphism $f:T\rightarrow T' $ is called a {\em weak equivalence}
\cite{Illusie} if for all $Y\in Sch /\cc$ and $t\in T(Y)$, the resulting
morphisms  $$
\pi _i(T|_{Sch /Y}, t)\rightarrow \pi _i(T'|_{Sch /Y}, f(t))
$$
are isomorphisms of sheaves on $Sch /Y$ (resp. $\pi _0(T)\rightarrow \pi
_0(T')$ is an isomorphism of sheaves of sets on $Sch /\cc $). The theory is
localized by this notion of equivalence, in other words $T$ and $T'$ are
thought of as equivalent if there is a weak equivalence between them. Jardine
constructs a closed model category which takes this into account \cite{Jardine}.
This leads, in particular, to the right notion of morphism, namely we only look
at morphisms whose target is a fibrant object; if necessary, a target object is
replaced by a weakly equivalent fibrant object. Without further mentionning
this, we make the convention that whenever we speak of morphisms between
$n$-stacks, the target object is made fibrant.

For a general site, one can have a connected stack $T$ (i.e. $\pi _0(T)=\ast$)
but where the global section space of $T$ is empty, or nonconnected. However,
in the present case we are working in the etale topology over an algebraically
closed field $\cc$. In this case we have the implication
$$
\pi _0(T)= \ast \; \Rightarrow \; \pi _0(T(Spec (\cc ))= \ast .
$$
Indeed, the etale coverings of $Spec (\cc )$ are trivial, so there is no change
over the object $Spec (\cc )$ when one passes from the presheaf $\pi _0^{\rm
pre}(T)$ to the associated sheaf.

If $T$ is connected, then, we can choose a basepoint $t\in T(Spec (\cc ))$ which
is unique up to homotopy, so the sheaf of groups $\pi _1(T,t)$ is uniquely
defined up to global conjugacy.  If $Y$ is any scheme and $t'\in T(Y)$ then
locally on $Y$, $t'$ is equivalent to $t|_Y$ so $\pi _1(T, t')$ is locally over
$Y$ equivalent to the restriction of $\pi _1(T,t)$ to $Y$. Thus in this case,
the fundamental group sheaf $\pi _1(T,t)$ over $Sch /\cc $ gives a relatively
accurate picture of the $1$-type of $T$.

In particular, it makes sense to require that $T$ be {\em $1$-connected}, that
is that $\pi _0(T)=\ast $ and $\pi _1(T,t)=\{ 1\} $ for the
basepoint $t\in T(Spec (\cc ))$. If $T$ is $1$-connected then for any scheme
$Y$ and $t'\in T(Y)$, $\pi _1(T,t')=\{ 1\}$. Furthermore, in this case
the fundamental group of $T(Spec (\cc ))$ is trivial (the cohomological
contributions from the higher homotopy vanish because etale cohomology of $Spec
(\cc )$ is trivial). Therefore the basepoint $t\in T(Spec (\cc ))$ is
well-defined up to unique homotopy.

We now describe the standard topological constructions which we shall use for
$n$-stacks. The first is the notion of {\em homotopy fiber product}. If
$A\rightarrow B \leftarrow C$ are morphisms of $n$-stacks then we obtain
the {\em homotopy fiber product}
$A\times _BC$
with a
diagram
$$
\begin{array}{ccc}
A\times _BC&\rightarrow & A\\
\downarrow && \downarrow \\
C & \rightarrow & B
\end{array}
$$
together with a homotopy of commutativity of the diagram. These data are
essentially well-defined (in the sense that they are well-defined up to
homotopy which is itself well-defined up to homotopy \ldots ). In the
simplicial presheaf theory, the homotopy fiber product is obtained by replacing
one of the two morphisms by a fibrant morphism and then taking the usual fiber
product. In the $n$-category theory, see \cite{limits}.

Suppose $f:A\rightarrow B$ is a morphism and suppose $b\in B(Spec (\cc ))$. We
can think of $b$ as a morphism $b: \ast \rightarrow B$ where $\ast$ denotes the
constant presheaf with values the $1$-point topological space. Define the {\em
fiber of $f$ over $b$} to be the homotopy fiber product
$$
Fib(f,b):= \ast \times _B A.
$$
If the base  $B$ is $0$-connected, then as mentionned above, the choice of
basepoint $b$ exists and is unique up to a global homotopy (i.e. a path in
$B(Spec (\cc ))$.
Thus we can denote by $Fib(f)$ the fiber over this $b$, bearing in mind that it
is defined up to the conjugation action of $\pi _1(B (Spec (\cc )))$.
If $B$ is $1$-connected the choice of basepoint $b$ is unique up to unique
homotopy, so $Fib(f)$ is well-defined up to homotopy. We call this the {\em
homotopy fiber of $f$}.

We say that
$$
A\rightarrow B \rightarrow C
$$
is a {\em fiber sequence} if $C$ is $1$-connected or if we are otherwise given
a basepoint $c\in C(Spec (\cc ))$, if we are given a homotopy between the
composition $A\rightarrow C$ and the constant  map at the basepoint $c$, and if
the map $A\rightarrow B$ together with this homotopy induce an equivalence
between $A$ and $Fib(B\rightarrow C, c)$.

A morphism $T\rightarrow R$ is said to be a {\em locally constant fibration
with fiber $F$} if for every scheme $Y\rightarrow R$, locally on $Y$ (in the
etale topology) we have $Y\times _RT \cong Y \times F$. In the usual
topological case with connected base, this is vacuous. In the case of
$n$-stacks over a site, we  still have that if $R$ is $0$-connected then any
morphism $T\rightarrow R$ is a locally constant fibration. However, we are
often interested in cases where $R$ is not connected (i.e. $\pi _0(R)$ is
some nontrivial sheaf of sets). In these cases, being locally constant
condition is a nontrivial additional condition.

The next general type of operation we discuss is {\em truncation}. If $T$ is an
$n$-stack and if $m\leq n$ then we obtain an $m$-stack $\tau_{\leq m}T$
(i.e. a simplicial  presheaf which is $m$-truncated, in other words has
homotopy group sheaves vanishing in degrees $>m$) together with a morphism of
$n$-stacks $$
T\rightarrow \tau_{\leq m}T
$$
which induces an isomorphism on homotopy group sheaves in degrees $i\leq m$.
The $m$-stack $\tau_{\leq m}T$ together with this morphism are essentially
well-defined. We can construct $\tau _{\leq m}T$ as the $m$-stack associated to
the presheaf of spaces
$$
(\tau _{\leq m}^{\rm pre}T)(Y):=
\tau _{\leq m}(T(Y))
$$
where the truncation on the left is just truncation of topological spaces
(also known as the coskeleton operation).

A first example of truncation is the sheaf of sets $\pi _0(T)= \tau _{\leq 0}T$.

If $T\rightarrow R$ is a morphism of $n$-stacks then there is a relative
(or ``fiberwise'')
version of the truncation denoted $\tau _{\leq m/R}(T)\rightarrow R$. This is
defined by the property that for any scheme $Y$ and morphism $Y\rightarrow R$,
$$
\tau _{\leq m/R}(T)\times _R Y = \tau _{\leq m}(T\times _RY).
$$

Using the operations of truncation and homotopy fiber products, we obtain the
{\em Postnikov tower}. If $T$ is an $n$-stack then we have morphisms
$$
T\rightarrow \ldots \rightarrow \tau _{\leq m}T
$$
$$
\rightarrow \tau _{\leq m-1}T \rightarrow \ldots \rightarrow \pi _0(T).
$$
In order to describe the stages in this tower of maps, we need a few more
notions.

Suppose $Y$ is a scheme and $L$ is a sheaf of groups over $Y$. Fix $m\leq n$ and
suppose $L$ is abelian if $m\geq 1$. Then we can construct the simplicial
presheaf $K^{\rm pre}(L,m)$ on $Sch/Y$ by the standard construction applied
to $L$; let $K(L,m)$ be the associated stack.  Note that
$K(L,m)$ has a chosen basepoint section (over $Y$) which we denote by $0$, and
$\pi _i(K(L,m),0)=0$ for $i\neq m$, and it is $=L$ for $i=m$. Furthermore these
properties characterize $K(L,m)$ essentially uniquely.

The $K(L,m)$ on the site $Sch /Y$  corresponds to an $n$-stack
on $Sch /\cc$ with morphism to $Y$ (where $Y$ is considered as a $0$-stack or
sheaf of sets), which we denote by $K(L/Y,m)\rightarrow Y$.

We can do the same construction relative to any $n$-stack but for this we need
to have a notion corresponding to sheaf of (abelian) groups.
If $A$ is an $n$-stack then a {\em local system of (abelian) groups)} on
$A$ is a morphism $L\rightarrow A$ with relative group  structure, which is
relatively $0$-truncated (i.e. for any scheme $Y$ and map $Y\rightarrow A$, the
homotopy fiber product $Y\times _AL$ is $0$-truncated). This is equivalent to
the data for every $Y$ of a local system of (abelian) groups $L_Y$ over $A(Y)$,
together with restriction morphisms $L_Y|_{A_{Y'}}\rightarrow L_{Y'}$
for $Y'\rightarrow Y$,
satisfying the obvious associativity condition.

{\em Caution:} if $X$ is a sheaf of sets represented by a scheme, then
a local system over $X$ (according to the above terminology) is the same
thing as
a sheaf of (abelian) groups over $X$. It doesn't have anything to do with the
notion of ``flat vector bundle'' over $X$.

If $L\rightarrow A$ is a local
system of abelian groups then we obtain a morphism
$$
K(L/A,n)\rightarrow A,
$$
whose homotopy fiber over any $a\in A(Y)$ is the Eilenberg-MacLane $n$-stack
$K(L|_Y,n)$ over $Y$. For $n=1$ we can make do with any local
system of groups not necessarily abelian.

There is a standard fibration sequence relative to $A$
$$
K(L/A,m)\rightarrow A \rightarrow K(L/A, m+1),
$$
in other words
$$
K(L/A,m)= A \times _{K(L/A, m+1)}A.
$$

Using this we obtain the usual description of the stages in the Postnikov
tower: if $m\geq 2$ then, setting $A:= \tau _{\leq m-1}T$ there is a local
system $L$ of abelian groups over $A$ and a section $ob:A\rightarrow K(L/A,
m+1)$
such that the morphism  in the Postnikov tower
$$
\tau _{\leq m}T \rightarrow \tau _{\leq m-1}T =A
$$
is equivalent to
$$
A\times _{K(L/A,m+1)} A\rightarrow A
$$
where the first morphism in the fiber product is $ob$ and the second is $0$.

The description of the first stage $\tau _{\leq 1}T\rightarrow \pi _0(T)$
is much more complicated and is basically the subject of Giraud's book
\cite{Giraud}.

Using the notion of local system we can define a relative version of the
homotopy group sheaves. If $T\rightarrow R$ is  a morphism of $n$-stacks and $s:
R\rightarrow T$ is a section then we obtain local systems of groups $\pi
_i(T/R, s)$ over $R$.

Suppose $X\rightarrow Z$ and $Y\rightarrow Z$ are morphisms of $n$-stacks. Then
we obtain a {\em relative internal $Hom$} which is an $n$-stack with morphism
to $S$,
$$
Hom (X/Z,Y/Z)\rightarrow S.
$$
It is defined by the universal property that maps $A\rightarrow Hom (X/Z,Y/Z)$
are the same (in an essentially well-defined way) as maps $X\times
_ZA\rightarrow Y$ over $Z$.  For existence, if the proof
isn't contained somewhere in the literature then one might have to apply the
techniques of \cite{limits}. If $Z=\ast$ then we get back to the usual internal
$Hom(X,Y)$.

Similarly if  $Y\rightarrow X \rightarrow Z$ then we obtain the {\em relative
section stack}
$$
\Gamma (X/Z, Y)
$$
which is defined to be the fiber product
$$
Hom (X/Z,Y/Z)\times _{Hom (Y/Z,Y/Z)}Z
$$
where the second map in the fiber product is that corresponding to the identity
of $Y$, and the first map is induced by $Y\rightarrow X$. Again if $Z=\ast$ we
denote this simply by $\Gamma (X,Y)$.

We now come to one of the main types of observations, namely the relationship
between the above objects and cohomology. See for example Thomason
\cite{Thomason} for much of this. If $A$ is an $n$-stack and $L$ a local system
of abelian groups over $A$ then we define
$$
H^i(A, L):= \pi _0\Gamma (A, K(L/A,i)).
$$
It is a sheaf of abelian groups on the site $Sch /\cc $. Similarly if $Z$ is an
$n$-stack, $p:A\rightarrow Z$ a morphism and $L$ a local system of abelian
groups
over $A$ then we define
$$
H^i(A/Z,L):= \tau _{\leq 0 /Z}\Gamma (A/Z, K(L/A,i))
$$
where $\tau _{\leq 0 /Z}$ is the relative version of the truncation operation
for $n$-stacks over $Z$. Note that $H^i(A/Z,L)$ is a local system of abelian
groups on $Z$. We can also denote it by $R^ip_{\ast} (L)$.

One has the result that the cohomology defined above coincides with sheaf
cohomology over simplicial objects (representing $A$ by a simplicial object in
the topos of $Sch /\cc$). See \cite{Thomason} or
\cite{flexible}. In particular the notation $R^ip_{\ast} (L)$ coincides with the
usual meaning (particularly when we are looking at $A$ and $Z$ which are
represented by schemes, for example).

We have the formulae
$$
\pi _i(\Gamma (A, K(L/A,m)), 0) =H^{m-i}(A,L)
$$
and (using the relative version of homotopy groups)
$$
\pi _i(\Gamma (A/Z, K(L/A,m))/Z, 0) =H^{m-i}(A/Z,L).
$$

The usual results concerning cohomology of topological spaces hold for
cohomology as defined above. In particular, we have cup-products,
corresponding to the following operations on Eilenberg-MacLane spaces. If
$L$, $L'$  and $L''$ are local systems of abelian groups over $A$ and if
$$
L\times _AL'\rightarrow  L''
$$
is a bilinear morphism (of relative abelian group objects) then we obtain
morphisms  $$
K(L/A, i)\times K(L'/A, j)\rightarrow K(L'' /A, i+j).
$$
These give cup-products in cohomology which are bilinear morphisms
$$
H^i(A, L)\times H^j(A,L')\rightarrow H^{i+j}(A, L'').
$$

We also have a K\"unneth formula. The case which we use in the present paper is
as follows. Suppose $X$ and $Y$ are $n$-stacks. Then
$$
H^m(X\times Y, \Oo )= \bigoplus _{i+j= m} H^i(X,  H^j(Y,\Oo )).
$$
If $H^j(Y,\Oo )$ are represented by finite dimensional vector spaces then we can
write the more usual formula
$$
H^m(X\times Y, \Oo )= \bigoplus _{i+j= m} H^i(X, \Oo )\otimes _{\Oo} H^j(Y,\Oo
).
$$
This extends to the relative case of morphisms $X\rightarrow S$ and
$Y\rightarrow S$ if these families are locally trivial over the etale topology
of $S$.

Finally, we have a Leray-Serre spectral sequence. See \cite{Thomason} for
one way
to set this up. If $f: X\rightarrow Y$ is a morphism of $n$-stacks and if $L$ is
a local system of abelian groups on $X$ then we obtain a ``complex'' $R^{\cdot}
f_{\ast}(L)$ on $Y$ and the cohomology of $X$ is the ``hypercohomology'' of this
complex. These terms are put in quotations because one should actually interpret
the notion of complex as being a fibration in spectra over $Y$ (the raw notion
of complex of local systems is not adapted to the higher homotopy involved if
$Y$ is not $0$-truncated and locally cohomologically trivial). In any case we
get the cohomology objects of the direct image, which are the relative
homotopy group sheaves of the fibration of spectra, denoted $R^if_{\ast}(L)$.
These are local systems over $Y$. We have the Leray-Serre spectral sequence
(cf Thomason \cite{Thomason})
$$
E^{i,j}_2 = H^i(Y, R^jf_{\ast}(L))\Rightarrow H^{i+j}(X, L).
$$

One way of approaching all of the above details is to replace any
$n$-stack i.e. simplicial presheaf (particularly those coming in as the base in
questions about cohomology) by a simplicial object whose stages are formal,
possibly infinite, disjoint unions of schemes (this technique was pointed out to
me by C. Teleman). Similarly, we can replace morphisms of $n$-stacks by
morphisms of such objects. Then questions about cohomology become just
questions about cohomology of simplicial schemes with obvious modifications
made to allow for the infinite number of components.

The $n$-stack $K(\Oo , n)$ has an infinite loop-space structure or
$E_{\infty}$-structure, in other words it has
an infinite delooping (the $m$-fold delooping is just $K(\Oo , m+n)$).
This structure is the homotopical analogue of an abelian group structure:
it contains a homotopy class of maps
$$
K(\Oo , n)\times K(\Oo , n)\rightarrow K(\Oo ,n)
$$
but also higher homotopies of associativity, commutativity etc.
We often think of $K(\Oo , n)$ as a homotopical group object and speak of
things such as ``principal bundles'' for it.

The infinite loop-space structure is
inherited by $Hom (\Ff , K(\Oo , n))$ for  any sheaf $\Ff$.

In the first version of the paper, it was stated that the infinite loop-space
structure provides a functorial decomposition into products of
Eilenberg-MacLane stacks.  This is not true in general. For one thing, such a
decomposition may not exist, and for another it is never completely functorial.
For existence, the obstruction is in the $Ext$ between the various homotopy
group sheaves. If such a decomposition exists, then it can be functorial up to
one homotopy, but this homotopy itself will not be uniquely determined up to a
second homotopy.

The statement which we actually need for the calculation in \S 8 below, does
work, and we state it as a proposition.

\begin{proposition}
\label{decomp}
(A)\, Suppose $(Y,y)$ is a basepointed $n$-stack with infinite loop structure.
Suppose that $\pi _i(Y,y)$ are represented by finite-dimensional vector spaces.
Then there exists an equivalence
$$
\varepsilon : Y\cong Y_0 \times \ldots \times Y_n
$$
where
$$
Y_i = K( \Oo ^{a_i}, i) = K(\pi _i(Y,y), i).
$$
(B)\, Suppose $f:Y\rightarrow Y'$ is a map of basepointed $n$-stacks both as in
(A), and suppose $f$ is compatible with the infinite loop structure. Let
$$
\varepsilon : Y\cong Y_0\times \ldots \times Y_n,\;\;\;
\varepsilon ': Y'\cong Y'_0\times \ldots \times Y'_n
$$
denote the maps given by (A) (chosen independantly of $f$). Then there is a
homotopy making the square
$$
\begin{array}{ccc}
Y&\cong &Y_0\times \ldots \times Y_n\\
\downarrow &&\downarrow \\
Y'&\cong &Y'_0\times \ldots \times Y'_n
\end{array}
$$
commute, where the vertical arrow on the right is a product of the morphisms
of Eilenberg-MacLane stacks
$Y_i\rightarrow Y'_i$ induced by
$f_{\ast}:\pi _i(Y,y)\rightarrow \pi _i(Y', y')$.
\end{proposition}
{\em Proof:}
The short way of saying this is that the $Ext^j(V,W)$ vanish for $j>0$ for
finite dimensional vector spaces, so the obstruction to splitting vanishes.
We give the following more concrete argument (which is in a certain sense just
repeating the argument which will be given in \ref{ext} below for the vanishing
of the $Ext^j$).

For $N>n$ we are given an $N-1$-connected pointed $n+N$-stack $(Z,z)$ with
$(Y,y)= \Omega ^N(Z,z)$ (this is the {\em ad hoc} definition of ``infinite loop
structure'' which we use). In part (B) the map $f$ comes from a map
$g:Z\rightarrow Z'$.  Thus for part (A) it suffices to obtain a decomposition
$$
Z\cong Z_0\times \ldots \times Z_n
$$
with
$$
Z_i = K(V_i, N+i),\;\;\; V_i := \pi _i(Y,y) = \pi _{N+i}(Z,z),
$$
and for part (B) it suffices
to obtain the homotopy of functoriality on the level of the morphism $g$.

Calculate the cohomology of $Z$ with coefficients in $\Oo$. Using Leray-Serre
spectral sequences for the stages in the Postnikov tower, and using the Breen
calculations \ref{bc} in view of the hypothesis that the $\pi _j$ are finite
dimensional vector spaces, we find that for $i\leq N+n$,
$$
H^{j}(Z, \Oo ) = Hom (\pi _j(Z,z), \Oo ) =V_{N-j}^{\ast}.
$$
Thus we have tautological classes
$$
\varepsilon _i \in H^{N+i}(Z, \pi _{N+i}(Z,z))= H^{N+i}(Z, V_i)
$$
which together provide us with a map
$$
\varepsilon = (\varepsilon _0,\ldots , \varepsilon _n) :
Z \rightarrow K(V_0,N)\times \ldots \times K(V_n, N+n).
$$
This map induces an isomorphism on homotopy groups (in degrees up to $N+n$).
Thus it provides the required splitting for (A).

For part (B) suppose we have a map
$$
g:Z= Z_0\times \ldots \times Z_n \rightarrow Z'_0\times \ldots \times Z'_n
$$
with $Z_i= K(V_i, N+i)$ and $Z'_i= K(W_i, N+i)$. Such a map corresponds, up to
homotopy, to a collection of classes in
$H^{N+i}(Z, W_i)$.  From the K\"unneth formula (which can be seen by a
collection
of Leray-Serre spectral sequences for the projections onto the factors) we have
$$
H^{N+i}(Z, W_i) = H^{N+i}(Z_i, W_i).
$$
The other factors vanish again using the Breen calculations \ref{bc} from the
fact that $V_j$ and $W_i$ are represented by finite dimensional vector spaces.
Our map $g$ is therefore homotopic to a map given by the classes
in $H^{N+i}(Z_i, W_i)$, i.e. a map compatible with the product decomposition.
\eop

The homotopy in part (B) is not unique: it can be changed by a map
$Y\rightarrow \Omega Y'$ in other words by a collection of morphisms
$Y_i\rightarrow Y'_{i+1}$.

\medskip

Suppose that $\Ff$ is an $n$-stack such that the $H^i(\Ff , \Oo )$ are
represented by finite dimensional vector spaces. Then we can apply the above
proposition to
$$
Y:=  Hom (\Ff , K(\Oo ,n)).
$$
The decomposition of $Hom (\Ff , K(\Oo ,n))$ into a
product of Eilenberg-MacLane spaces
$$
Hom (\Ff , K(\Oo ,n)) = \prod _i K(H^{n-i}(\Ff , \Oo ), i),
$$
is related to the K\"unneth formula.
A morphism $Z\rightarrow Hom (\Ff , K(\Oo , n))$ corresponds (by the definition
of internal $Hom$) to a morphism
$Z\times \Ff \rightarrow K(\Oo ,n)$, in other words to a class $f\in H^n(Z\times
\Ff , \Oo )$. By the above product structure this class decomposes into a
collection of classes $f_i\in H^i(Z, H^{n-i}(\Ff , \Oo ))$. The $f_i$ are the
K\"unneth components of $f$.

\begin{center}
$\ast$ \hspace*{2cm}$\ast$ \hspace*{2cm}$\ast$
\end{center}

For the remainder of the paper, we look at $n$-stacks of groupoids on $Sch /\cc$
and unless specified otherwise, the reader may fix any $n\geq 3$ (for the
calculation it suffices to take $n=3$.)

\numero{Connected very presentable shape}

We isolate some special $n$-stacks $T$ and then use them to measure the
``shape''
of an arbitrary $n$-stack $\Ff$. In other words, take a sub-$n+1$-category
${\cal P} \subset nSTACK$ of the $n+1$-category of $n$-stacks (of groupoids,
say), and look at the $n+1$-functor
$$
{\cal P}  \rightarrow  nSTACK
$$
$$
T\mapsto Hom (\Ff , T).
$$
We call this $n+1$-functor the {\em shape of $\Ff$ as measured by ${\cal P}$}.
There are many possible ways to choose ${\cal P}$. Some reasonable parameters
are to require that ${\cal P} \subset nSTACK$ be a full sub-$n+1$-category
(in other words that we make no limitation on the morphisms of ${\cal P}$);
and that the condition $T\in {\cal P}$ should be measured only by looking at
the homotopy group sheaves $\pi _i(T,t)$. For our present purposes we start by
requiring that $T$ be connected, i.e. $\pi _0(T)=\ast$. Among other things, this
insures that the isomorphism classes of the higher homotopy group sheaves $\pi
_i(T,t)$ be well defined.

Before describing our choice of conditions for the $\pi _i(T,t)$ we take note
of the following: any $T\in {\cal P}$ will decompose in a Postnikov tower
whose stages are $K(\pi _i , i)$.  Morphisms $\Ff \rightarrow K(\pi _i,i)$
are classified by $H^i(\Ff ,\pi _i)$ and more generally, one has obstruction
theory for classifying the morphisms $\Ff \rightarrow T$ going up in the
Postnikov tower; the obstruction classes are in  $H^{i+1}(\Ff ,\pi _i)$.
Thus, one should choose the class of possible $\pi _i$ to be a class of sheaves
such that, for the $\Ff$ we are interested in, the cohomology $H^j(\Ff ,\pi _i)$
has reasonable properties.

For the topic of Dolbeault cohomology, we already know how to take nonabelian
$H^1$ with coefficients in an affine algebraic group \cite{Moduli} \cite{hbls},
and we know how to take higher Dolbeault cohomology with coefficients in $\cc$
or more generally in a finite-dimensional complex vector space. This suggests
that our condtions should be that $\pi _1$ be an affine algebraic group, and
$\pi _i$ be represented by finite-dimensional vector spaces for $i\geq 2$.

Recall from \cite{RelativeLie} and \cite{GeometricN} that a {\em connected very
presentable $n$-stack $T$} is an $n$-stack of groupoids $T$ on $Sch /\cc$
subject
to the following conditions:
\newline
{\bf (connectedness):} $\pi _0(T) = \ast$ as a sheaf of sets on $Sch /\cc$;
\newline
{\bf (very presentability):} if $t\in T(Spec (\cc ))$ is a basepoint
(which we assume exists) then $\pi _i(T,t)$ are representable by group schemes
of finite type over $Spec (\cc )$, which are required to be affine for $i=1$
and vector spaces (i.e. affine unipotent abelian) for $i\geq 2$.

(This is the ``very presentability'' condition of \cite{RelativeLie} under the
additional hypothesis of connectedness; in the non-connected case the
definition is more complicated and that is basically the subject of the paper
\cite{RelativeLie}.)

We now choose ${\cal P}$ for our shape theory to be the $n+1$-category of
connected very presentable $n$-stacks of groupoids. Suppose $\Ff$ is an
$n$-stack on $Sch/\cc$. The{\em shape of $\Ff$} is defined as the
$n+1$-functor from the $n+1$-category of connected very presentable $n$-stacks
$T$, to the $n+1$-category $nSTACK$ of all $n$-stacks, given by the formula
$$
Shape (\Ff )(T):= Hom (\Ff , T).
$$
This contains all information about $\Ff$ which one can extract by looking at
$G$-torsors over $\Ff$ and cohomology of
associated vector bundles.

In many cases, $Hom(\Ff , T)$ will be a geometric or locally geometric
$n$-stack \cite{GeometricN}. For example in the case $\Ff = X_{Dol}$ we look at
below, $Hom(X_{Dol} , T)$ will be locally geometric. In these cases the shape
of $\Ff$ may be considered as an $n+1$-functor from connected very presentable
$n$-stacks to the $n+1$-category of (locally) geometric $n$-stacks, sitting
inside $nSTACK$.

\subnumero{Examples}

We explain how to understand the structure of connected very presentable $T$.
First is the simply connected case. Here $T$ is given by a Postnikov tower
where the stages are of the form $K(\Oo ^a, m)$.  The only question is how they
are put together. The fibration
$$
K(\Oo ^a, m)\rightarrow \tau _{\leq m}T \rightarrow \tau _{\leq m-1}T
$$
is classified by a map
$$
\tau _{\leq m-1}T\rightarrow K(\Oo ^a, m+1),
$$
in other words $\tau _{\leq m}T$ is the pullback by this map of the standard
fibration
$$
K(\Oo ^a, m)\rightarrow \ast \rightarrow K(\Oo ^a, m+1).
$$
We can write
$$
\tau _{\leq m}T= \tau _{\leq m-1}T\times _{K(\Oo ^a, m+1)}\;\; \ast .
$$
The classifying map is a class in $H^{m+1}(\tau _{\leq m-1}T, \Oo ^a)$.
In turn, this cohomology can be ``calculated'' by the Leray-Serre spectral
sequence applied to the previous part of the Postnikov tower for
$\tau _{\leq m-1}T$. The basic pieces that we need to know are the cohomology
of the Eilenberg-MacLane spaces. These are given by the {\em Breen
calculations} \cite{Breen1} \cite{Breen2}, which we recall in Appendix I
(giving a relative version). For the present discussion the answer is that
$H^{\ast} (K(\Oo ^a, m), \Oo )$
is a graded-symmetric algebra on $\Oo^a$ in degree $m$. Note that this answer
is the same as the classical answer for rational cohomology of rational
Eilenberg-MacLane spaces $H^{\ast} (K(\qq ^a, m),\qq )$.

As we shall explain below and also in Appendix II, if $Y$ is a finite
simply-connected $CW$-complex then we obtain a $1$-connected very presentable
$T= Y\otimes \cc$ whose homotopy group sheaves are $\pi _i(Y,y)\otimes _{\qq}
\Oo$.

We now look at connected but not simply connected very presentable $T$.
Note that
since $T$ is connected, we can choose a basepoint $t\in T(Spec (\cc ))$. Let
$G:= \pi _1(T,t)$; by hypothesis it is an affine algebraic group scheme over
$\cc$. We have a fiber sequence
$$
T'\rightarrow T \rightarrow K(G,1)
$$
where $T'$ is simply connected very presentable. The homotopy group sheaves
$\pi _i(T', t)= \pi _i(T,t)$ are vector spaces $\Oo ^a$ for $i\geq 2$; but
notice also that $G$ acts on these vector spaces. The action is an action of
sheaves on the site $Sch /\cc $ so it is automatically algebraic; we can write
$\pi _i(T', t)= V^i$ with $V^i$ a linear representation of $G$.

The same Postnikov tower discussion as above, works here. The only difference
is that in calculating the cohomology of the $\tau _{\leq m}T$ we may have
coefficients which are linear representations of $G$, and at the end we get
down to a step where we have to calculate $H^i(K(G,1), V)$. If $G$ is
reductive, this ``algebraic cohomology'' vanishes for $i\geq 1$, whereas if $G$
is unipotent then it is equal to the Lie algebra cohomology.

A simple example of non-simply connected $T$ may be obtained as follows. Start
with a linear algebraic group $G$ with a linear representation $V$,
corresponding to a local system $\underline{V}\rightarrow K(G,1)$.
Use the notational shorthand
$$
K(V/G; n):= K(\underline{V}/K(G,1), n).
$$
We have a fibration sequence
$$
K(V, n)\rightarrow K(V/G,n) \rightarrow K(G,1).
$$
Maps $\Ff \rightarrow T$ correspond
to pairs $(E, \eta )$ where $E$ is a $G$-torsor over $\Ff$ and
\newline
$\eta \in H^n(\Ff,
E\times ^GV)$.

An advantage of the nonabelian cohomological formulation of the
secondary Kodaira-Spencer classes we define is that the definition works for
cohomology with coefficients in a connected very presentable $T$, even
non-simply connected. However, for the calculation we will do, we look at a
particular simply-connected $T$ (the ``complexified $2$-sphere''
$S^2\otimes \cc$).

\subnumero{Representability}

Under certain circumstances, basically when the shape of $\Ff$ is simply
connected and has reasonable cohomology sheaves, then $Shape (\Ff )$ is {\em
representable}. By this we mean that there is a morphism $\Ff \rightarrow
\Sigma$ from $\Ff$ to a very presentable $n$-stack $\Sigma$ such that for any
other very presentable $n$-stack $T$ we have
$$
Hom (\Sigma , T)\stackrel{\cong}{\rightarrow} Hom (\Ff , T).
$$
For example we have the following precise statement.

\begin{theorem}
\mylabel{representable0}
Suppose $\Ff$ is an $n$-stack on $Sch /\cc $ such that for any affine algebraic
group $G$,
$$
K(G,1)\stackrel{\cong}{\rightarrow}Hom (\Ff , K(G,1)).
$$
Suppose that the $H^i(\Ff , \Oo )$ are representable by finite dimensional
vector spaces. Then there is a morphism $\Ff \rightarrow \Sigma$ to a
$1$-connected very presentable $n$-stack $\Sigma$, such that for any connected
very presentable $n$-stack $T$ we have
$$
Hom (\Sigma , T)\stackrel{\cong}{\rightarrow} Hom (\Ff , T).
$$
\end{theorem}
The proof will be given in Appendix II.

On the other hand, for $n$-stacks $\Ff$ whose shape is not $1$-connected, the
shape will not in general be representable. For example, suppose $W$ is a
finite $CW$ complex (considered as a constant $n$-stack) such that $\pi
_1(W)=\Gamma := \zz $ and such that some $\pi _i(W)\otimes \qq$ is a $\qq
[\Gamma
]$-module which is not completely torsion (hence infinite-dimensional over
$\qq$). Then $Shape (W)$ is not representable. Indeed, the infinite
dimensionality of $\pi _i(W)\otimes \qq$ is seen by the shape, since all
irreducible representations of $\Gamma$ are finite ($1$-) dimensional.

\numero{Nonabelian Dolbeault cohomology}

Suppose $X$ is a smooth quasiprojective variety. Recall that one defines the
{\em Dolbeault cohomology} of $X$ as the hypercohomology of the trivial complex
$\Omega ^{\cdot}_X$ with differential equal to $0$:
$$
H^i_{Dol}(X):= {\bf H}^i(X, \Omega ^0_X \stackrel{0}{\rightarrow} \ldots )
= \bigoplus _{p+q=i} H^q(X, \Omega ^q_X)
$$
(generally speaking this is only motivated by topology when $X$ is projective;
but we make the notation for quasiprojective $X$ too, for use in Mayer-Vietoris
arguments).

A nonabelian version for $H^1$ may be defined by setting $H^1_{Dol}(X, G)$
equal to the moduli stack of Higgs principal $G$-bundles (\cite{Hitchin0}
\cite{Hitchin}
\cite{hbls} \cite{Moduli}) $(P,\theta )$.
One can impose semistability conditions
and vanishing of Chern classes to get a version more closely related to
topology,
but we don't need that for the present algebraic discussion.

If $(P,\theta )$ is a principal Higgs bundle and if $V$ is a representation of
$G$ then we obtain an associated Higgs bundle $(E,\theta )$. Recall that we
define the {\em Dolbeault cohomology of $(E,\theta )$} as the hypercohomology of
the {\em Dolbeault complex} (cf \cite{hbls})
$$
\ldots \stackrel{\theta}{\rightarrow} E\otimes _{\Oo _X}\Omega ^i_X
\stackrel{\theta}{\rightarrow} \ldots ,
$$
$$
H^i_{Dol}(X, (E,\theta )):= {\bf H}^i(X,
(E\otimes _{\Oo_X}\Omega ^{\cdot}_X,\theta )).
$$

We present a way of unifying these definitions into a notion of {\em nonabelian
Dolbeault cohomology}. For nonabelian $H^1$ the present interpretation was
explained in \cite{SantaCruz}.

Let $\widehat{TX}$ denote the formal completion of
$TX$ along the zero section. Considered as a presheaf on $Sch /\cc $ it
associates to any $\cc$-scheme $Y$, the set of maps $Y\rightarrow TX$ which map
the underlying reduced subscheme $Y^{\rm red}$ to the zero-section $X\subset
TX$.

Define the
$1$-stack $X_{Dol}$ to be the relative $K(\widehat{TX}/X, 1)$, i.e. the relative
classifying stack for the group scheme $\widehat{TX} \rightarrow X$.

A variant which is technically easier to work with is
$$
X_{UDol} := K(TX/X, 1).
$$
The $U$ in the notation stands for ``unipotent'': as we shall see below
(Proposition \ref{calcDol}), a morphism $X_{UDol} \rightarrow K(G,1)$
corresponds
to a principal Higgs bundle with structure group $G$, over $X$, such that the
Higgs field is a section of unipotent elements of the Lie algebra of $G$.
However, for morphisms to simply connected $T$, we can safely replace $X_{Dol}$
by $X_{UDol}$ and for the purposes of the present paper this is what we shall
do.

If $T$ is an $n$-stack then we define the {\em
nonabelian Dolbeault cohomology of $X$ with coefficients in $T$} to be the
$n$-stack (of $n$-groupoids)
$$
Hom (X_{Dol}, T).
$$
Of course this includes the classical abelian case, when we take $T= K(\Oo , n)$
where $\Oo$ denotes the structural sheaf, equal to ${\bf G}_a$, represented by
the affine line. It also includes the somewhat classical case of nonabelian
Dolbeault $H^1$ with coefficients in a group scheme $G$.

More generally we will be most interested in the case where $T$ is a connected
very presentable $n$-stack.
The calculation of the present paper involves a $1$-connected $T$, so it doesn't
refer to the case of Higgs principal $G$-bundles or Dolbeault cohomology with
coefficients in Higgs bundles. These aspects are only presented to show the
unified character of the definition.

\begin{proposition}
\mylabel{calcDol}
We have that
$$
\pi _0 (Hom (X_{Dol}, K(\Oo , n))) = H^n_{Dol}(X, \cc )
$$
is the usual Dolbeault cohomology of $X$. The same holds for
$X_{UDol}$. If $G$ is an affine algebraic group scheme then
$$
Hom (X_{Dol} , K(G,1))= {\cal M}_{Dol}(X, G)
$$
is the moduli $1$-stack of principal Higgs bundles with structure group $G$.
On the other hand,
$$
Hom (X_{UDol} , K(G,1))
$$
is the moduli $1$-stack of principal Higgs bundles $(P, \theta )$ with structure
group $G$, such that for every $x\in X$ the element $\theta _x \in ad(P)_x\cong
\gerg $ is a unipotent element of the  Lie algebra $\gerg$ of $G$ (this condition
is of course independant of the isomorphism $ad(P)_x\cong
\gerg $ chosen).
\end{proposition}
{\em Proof:}
Let $X^{\rm fc}_{Dol}$ denote the formal category defined in \cite{SantaCruz}
which gives the $1$-stack $X_{Dol}$. We have
$$
Ob\, X^{\rm fc}_{Dol} = X,
$$
and
$$
Mor (X^{\rm fc}_{Dol}) = \widehat{TX} \rightarrow X \hookrightarrow X\times X
$$
(the morphism object lies over the diagonal in $X\times X$). The composition of
morphisms is just addition in $\widehat{TX}$ (which is a formal group scheme
over $X$).  This formal groupoid gives in an obvious way a presheaf of
groupoids on $Sch /\cc $, whose associated stack is $X_{Dol}$. A morphism
$X_{Dol}\rightarrow K(G,1)$ is the same thing as a $G$-torsor over
$X_{Dol}^{\rm fc}$, which in turn is the same thing as a principal $G$-bundle
$P$ over $X$ together with action of the formal group scheme $\widehat{TX}$.
The action may be interpreted as a morphism of sheaves of groups over $X$,
$$
\widehat{TX} \rightarrow Ad (P).
$$
Since the domain is formal, this is the same thing as a morphism of Lie
algebras over $X$,
$$
\theta : TX \rightarrow ad(P).
$$
This proves the second statement. For the third statement, the same proof works
but with $\widehat{TX}$ replaced by $TX$. Note that $TX$ is a unipotent group
scheme over $X$, so
$$
TX\rightarrow Ad(P)
$$
corresponds to a morphism of Lie algebras
$$
\theta : TX \rightarrow ad(P)
$$
with image in the unipotent elements of $ad(P)$.

For the first statement, one way to proceed is to notice that the cohomology of
$X_{Dol}$ with coefficients in $\Oo$ is the same as the cohomology of the
formal category $X^{\rm fc}_{Dol}$ as considered by Berthelot \cite{Berthelot}
and Illusie
\cite{Illusie}. From those references, one gets a generalized de Rham complex
calculating the cohomology, which in our case is seen to be exactly the
Dolbeault complex.

We also need to prove that the morphism $X_{Dol}\rightarrow X_{UDol}$ induces an
isomorphism on cohomology. For this it suffices to look locally over $X$, so we
can assume that $TX$ is trivial. Thus it suffices to prove that if $V$ is a
vector space of dimension $n$ then the morphism
$$
K(\widehat{V},1)\rightarrow K(V,1)
$$
induces an isomorphism of cohomology. This morphism has homotopy fiber the
sheaf $V_{DR}$ defined by $V_{DR}(Y)=V(Y^{\rm red})$ (cf \cite{kobe}). The
cohomology of $V_{DR}$ with coefficients in $\Oo$ is the algebraic de Rham
cohomology of $V$ (\cite{kobe} Theorem 6.2) which is trivial because $V$ is an
affine space.

Since the statement of the first part for $X_{UDol}$ is what we actually use,
we indicate a somewhat more elementary proof. Let $V$ be a vector space of
dimension $n=dim(X)$. We have a natural (split) extension
$$
1\rightarrow V \rightarrow G \rightarrow GL(V) \rightarrow 1
$$
which gives a fiber sequence of $1$-stacks
$$
K(V,1)\rightarrow K(G,1)\stackrel{p}{\rightarrow} K(GL(V), 1).
$$
Let $C^{\cdot} = R^{\cdot} p_{\ast}(\Oo )$.  Using the Breen calculations
\cite{Breen2} \cite{kobe},
which we recall in Theorem \ref{bc} in Appendix I below, it is easy to see that
$$
H^i(C^{\cdot}) = \bigwedge ^i(V^{\ast}).
$$
However, the fact that $GL(V)$ is a reductive group implies that all of the
algebraic group cohomology, in other words the cohomology of $K(GL(V), 1)$
with coefficients in any local system associated to a representation of
$GL(V)$, vanishes except in degree $0$. Therefore the invariants of the complex
$C^{\cdot}$ vanish so
$$
C^{\cdot} \sim \bigoplus _{i} \bigwedge ^i(V^{\ast}).
$$
Now notice that the vector bundle $TX\rightarrow X$ corresponds to a map
$X\rightarrow K(GL(V), 1)$ and
$$
X_{UDol} = K(G,1)\times _{K(GL(V),1)} X.
$$
Let $q: X_{UDol}\rightarrow X$ denote the projection; then
$R^{\cdot} q_{\ast}\Oo $ is the pullback of $C^{\cdot}$, hence it splits as
$$
R^{\cdot} q_{\ast}\Oo = \bigoplus _i\Omega ^{i}_X.
$$
The cohomology of $X_{UDol}$ with coefficients in $\Oo$ is equal to the
hypercohomology of $X$ with coefficients in this complex, which is the Dolbeault
cohomology.
\eop

{\bf Remark:} We have the following which are sometimes useful:
$$
\pi _i (Hom (X_{Dol}, K(\Oo , n))) = H^{n-i}_{Dol}(X, \cc ).
$$

\begin{proposition}
\mylabel{calcDol2}
Suppose $(P,\theta )$ is a Higgs principal $G$-bundle on $X$ corresponding to a
map $X_{Dol}\rightarrow K(G,1)$. Suppose $V$ is a representation of $G$ and let
$(E,\theta )$ be the associated Higgs bundle. Then the cohomology of $X_{Dol}$
with coefficients in the local system $(P,\theta )^{\ast}(V)$ is naturally
isomorphic to the Dolbeault cohomology $H^n_{Dol}(X, (E,\theta ))$. The same
works for a family of principal bundles parametrized by a base scheme $S$.
\end{proposition}
{\em Proof:}
As before, we interpret $(P,\theta )$ as a $G$-torsor over the formal category
$X_{Dol}^{\rm fc}$ which defines the stack $X_{Dol}$. Associated to the
representation $V$ we get a local system over the formal category, and its
cohomology (which is that mentionned in the statement of the proposition) is
calculated by a de Rham complex (cf \cite{Berthelot} \cite{Illusie}). This de
Rham complex is exactly the Dolbeault complex for $(E,\theta )$.
\eop

\begin{corollary}
\mylabel{UdolDol}
Suppose $T$ is a connected very presentable $n$-stack such that $\pi _1(T)$ is
a unipotent group scheme over $Spec (\cc )$. Then the morphism
$$
Hom (X_{UDol} , T)\rightarrow Hom (X_{Dol} , T)
$$
is an equivalence.
\end{corollary}
{\em Proof:}
This follows immediately from the previous propositions using a Postnikov tower
for $T$, and noting that if $G$ is a unipotent group scheme then
all principal $G$-Higgs bundles satisfy the condition that the Higgs field be a
section of unipotent elements.
\eop

\subnumero{Cohomology classes of subschemes}

A particular detail of (usual) Dolbeault cohomology which enters into our
calculations is the cohomology class of a subscheme. Suppose $X$ is a smooth
projective variety and suppose that $Z\subset X$ is a subscheme of codimension
$d$. Then we obtain the class
$$
[Z] \in H^d(X, \Omega ^d_X)\subset H^{2d}_{Dol}(X).
$$
In fact this class has a canonical trivialization over $U:=X-Z$, in other words
we are actually given a lifting to
$$
[Z] \in H^d((X,U), \Omega ^d_X)\subset H^{2d}((X_{Dol}, U_{Dol}),\Oo ).
$$

Suppose for example $d=1$.
Take an open covering $X=\bigcup U_i$ with $U_0=U=X-Z$, and the
remaining $U_i$ affine with defining equations $z_i\in \Oo (U_i)$ having zeros
of order one along $Z$ and nonvanishing elsewhere. Put
$$
g_{ij}:= \frac{dz_i}{z_i} - \frac{dz_j}{z_j}.
$$
This gives a $1$-cocycle which determines the cohomology class $[Z]$. Since
$U_0=X-Z$ is a part of the covering, it gives a class in the cohomology of the
pair $(X,X-Z)$.

The case of higher codimension is treated similarly. The only thing we need to
know is that if $\lambda , \lambda ' \in \Oo
(X)$ are regular functions such that $\lambda |_Z= \lambda '|_Z$, then $\lambda
[Z] = \lambda '[Z]$. In particular if $P$ is a point then $\lambda [P]=\lambda
(P)[P]$.

\numero{Group actions}

Suppose $W$ is a sheaf of groups on $Sch /\cc $, and suppose $R$ is an
$n$-stack on $Sch /\cc $. Then an {\em action of $W$ on $R$} is a morphism
$$
\rho : A \rightarrow K(W, 1)
$$
together with an identification of the fiber $\rho ^{-1}(0)$ (by which
we mean the homotopy fiber product $\{ 0\} \times _{K(W, 1)} A$) with $R$,
$$
R\cong \rho ^{-1}(0).
$$
To put this more briefly, an action of $W$ on $R$ is a fiber sequence
$$
R\rightarrow A \rightarrow K(W,1).
$$

If $T$ is another $n$-stack, and if $\rho$ is an action of $W$ on $R$
then we obtain an associated action $Hom (\rho , T)$ of $W$ on $Hom (R,T)$.
It is given by the relative $Hom$,
$$
Hom (A/K(W,1), T)\rightarrow K(W,1),
$$
whose fiber is naturally identified with $Hom(R,T)$.

{\bf Example 1:} If $W$ acts on a sheaf of sets $R$ in the usual sense,
then this
can be interpreted as an action in the above sense. The stack $A$
(which in this case is a $1$-stack) is the
stack-theoretical quotient $R/W$ with its canonical principal $W$-bundle
$R\rightarrow R/W$ which corresponds to a morphism $R/W \rightarrow K(W,1)$.

{\bf Example 2:} The group $\Gm$ acts on $\widehat{TX}$ (resp. $TX$) over $X$ by
scalar multiplication. Therefore $\Gm$ acts on $X_{Dol}$ (resp. $X_{UDol}$)
and on $Hom (X_{Dol}, T)$ (resp. $Hom(X_{UDol}, T)$). In the case $T=K(G,1)$
this is the usual action of $\Gm$ on the moduli stack of principal Higgs bundles
\cite{Hitchin} \cite{hbls} \cite{Moduli}. In the case $T=K(\Oo , n)$ this
action gives the decomposition of $H^i_{Dol}(X)$ into pieces $H^q(X, \Omega
^p_X)$.

{\bf Example 3:} A different example, more closely related to what we are
interested in, is the following. Suppose
$$
1\rightarrow V \rightarrow E \rightarrow W \rightarrow 1
$$
is an exact sequence of sheaves of groups on $Sch /\cc$. This gives a fibration
sequence
$$
K(V,1)\rightarrow K(E,1)\rightarrow K(W,1),
$$
hence by definition it is an action of $W$ on $K(V,1)$.

More generally if $X$ is a sheaf of sets and if $V\rightarrow X$
is a sheaf of groups over $X$, then we can look at the relative
Eilenberg-MacLane stack  $K(V/X, 1)\rightarrow X$.
Suppose
$$
1\rightarrow V \rightarrow E \rightarrow p^{\ast}(W)\rightarrow 1
$$
is an exact sequence of sheaves of groups on $X$ (where $p:X\rightarrow \ast$
denotes the projection). Then we obtain  a fibration sequence
$$
K(V/X, 1)\rightarrow K(E/X, 1)\rightarrow K(W,1),
$$
the latter map being the composition of the map induced by
the second map in the exact sequence, with the projection
$$
K(p^{\ast}W /X,1)= K(W\times X/X, 1)= K(W, 1)\times X \rightarrow K(W,1).
$$
Thus, again by definition, our exact sequence corresponds to an action of $W$
on $K(V/X,1)$ (lying over the trivial action of $W$ on $X$).

\subnumero{Secondary classes}

We first discuss classifying spaces.
If $R$ is an $n$-stack then $Aut (R)$ is an $n$-stack with ``group'' structure,
more precisely with the structure of a loop space.
To make this statement precise, we construct a pointed $n+1$-stack
$B\, Aut(R)$ with basepoint denoted $0$, with an equivalence
$$
\Omega (B\, Aut(R), 0) \cong Aut (R).
$$
Construct the presheaf of spaces $B^{\rm pre}\, Aut(R)$ as the realization of a
simplicial $n$-stack $B_{\cdot}Aut(R)$ defined as follows: let
$\overline{I}^{(i)}$ be the $1$-category with $i+1$
isomorphic objects $0,\ldots , i$ (which we shall call ``vertices''); then
set $B_iAut(R)$ equal to the
$n$-stack of morphisms $\overline{I}^{(i)}\rightarrow n\underline{STACK}$
sending
the vertices to the object $R\in nSTACK$. Note that the
component $B_iAut(R)$ is homotopic to $Aut(R)\times \ldots \times Aut(R)$.
This simplicial $n$-stack may be interpreted as a presheaf of $n+1$-categories.
Note that $Aut(R)$ is defined
to be the $n$-stack of morphisms
$\overline{I}\rightarrow n\underline{STACK}$ sending $0,1$ to $R$, and it is a
stack of $n$-groupoids.  It may safely be confused with a presheaf of
$n$-truncated spaces, and  $B_{\cdot}\, Aut(R)$ becomes a presheaf of
simplicial spaces. The component simplicial spaces satisfy Segal's condition (cf
\cite{Segal} \cite{Tamsamani} \cite{nCAT} \cite{effective}) so if we set $B^{\rm
pre}\, Aut(R)$ equal to the  realization into a presheaf of spaces then Segal's
Proposition 1.5 \cite{Segal} implies that the natural map
$$
Aut(R) = B_1Aut(R) \rightarrow \Omega (B^{\rm pre}\, Aut (R), 0)
$$
is an equivalence (object-by-object).
Finally, let $B\, Aut(R)$ be the $n+1$-stack associated to $B^{\rm pre}Aut(R)$.

If $R$ is an $n$-stack of groupoids,
the locally constant fibrations with fiber $R$
$$
R \rightarrow E \rightarrow S
$$
are classified exactly by maps $S\rightarrow B\, Aut(R)$. In other words, given
$S$ and $R$, the $n+1$-category of such fiber sequences is equivalent to the
$n+1$-category $Hom (S, B\, Aut(R))$.

In particular, an action of a sheaf of groups $W$ on an $n$-stack of groupoids
$R$ is the same thing as a morphism
$$
f:K(W,1)\rightarrow B\, Aut (R).
$$
This leads to the notion of {\em characteristic classes} for the action:
if $\Gg$ is a sheaf of groups and
if $c\in H^i(B\, Aut (R), \Gg )$ then for any action we can pull back $c$ to
obtain a class in $H^i(K(W,1), \Gg )$.

One can also obtain {\em secondary invariants}, which only become defined when
some primary invariants vanish. For example, suppose
$c\in H^i(B\, Aut (R), \Gg )$ is a cohomology class corresponding to a map
$$
c: B\, Aut(R) \rightarrow K(\Gg , i).
$$
Suppose that $f: K(W, 1)\rightarrow B\, Aut(R)$ is an action of $W$ on $R$,
such that $c\circ f : K(W, 1)\rightarrow K(\Gg, i)$ is trivial, i.e. homotopic
to the constant map at the basepoint.  Let $Fib(c)$ denote the homotopy fiber
of the morphism $c$ over $0$. A choice of trivialization $\psi$ of $c\circ f$
gives rise to a morphism
$$
f_{\psi} : K(W, 1)\rightarrow Fib(c).
$$
The choices of trivialization $\psi$ are, up to homotopy, classified by
$H^{i-1}(K(W,1), \Gg )$.  If $c'$ is a class in $H^j(Fib(c), \Gg ')$
then the composition $c'\circ f_{\psi}$ gives a class in $H^j(K(W, 1), \Gg ')$.
This is (a typical example of) a secondary characteristic class for the action.

We don't get any further into the general theory of secondary classes, because
insofar as the formalism is concerned, there is no difference between the
present
case of $n$-stacks and the classical topological version of the theory.

Instead, we give a more explicit description of the secondary classes that we
are interested in.  These concern only classes of degrees $1$ and $2$ on $B\,
Aut (R)$, so a more concrete discussion is possible.
The first step is to replace $R$ by its truncation down to a $1$-stack,
denoted $\tau _{\leq 1}R$.  We fix a point $r\in R (Spec (\cc ))$ (and denote
also by $r$ its image in $\tau _{\leq 1}(R)$).  Note that the ``group''
$Aut(R)$
acts on $\tau _{\leq 1}R$.  The first of the primary classes comes from a
morphism
$$
Aut(R)\rightarrow Aut(\pi _0R).
$$
On the left is an actual sheaf of groups $\Gg$, and this map corresponds
to a class in $H^1(B\, Aut(R), \Gg )$. Suppose that we have an action by $W$
such that  this class vanishes. This implies that $W$ fixes the point $r\in \pi
_0(R)(Spec (\cc ))$.  (Vanishing of this class actually implies that $W$ fixes
every point but for what follows we only need that it fixes $r$).  Let
$R^r$ denote the component of $R$ containing $r$, more precisely it is the
fiber product
$$
R^r:= \ast \times _{\pi _0(R)} R.
$$
Under our preliminary vanishing hypothesis, we obtain an action of $W$ on $R^r$.

Now $\tau _{\leq 1}(R^r)= K(G, 1)$ where $G= \pi _1(R,r)$ (the point $r$ is a
basepoint defined over $Spec (\cc )$ giving this trivialization of the $1$-gerb
$\tau _{\leq 1}(R^r)$---the basepoint $r$ corresponds to $0\in K(G,1)$).

For any group $G$ we have
a fibration sequence
$$
K(Center(G), 2)\rightarrow
B\, Aut (K(G,1)) \rightarrow K(Out(G), 1).
$$
In our case above, $W$ acts on $R^r$, hence on $\tau _{\leq 1}R^r=K(G,1)$.
The map
$$
K(W,1)\rightarrow B\, Aut(\tau _{\leq 1}R^r)= B\, Aut (K(G,1))
$$
thus projects first of all to a morphism of groups
$$
W\rightarrow Out \left( \pi _1(R, r )\right) .
$$
This is in a certain sense again a primary invariant. Suppose that this
invariant vanishes. Then, given that
$$
K(W,1)\rightarrow K(Out(G), 1)
$$
is a pointed map, there is a canonical homotopy of this map to the constant map
at the basepoint, so we canonically can identify our map
$$
K(W,1)\rightarrow B\, Aut(\tau _{\leq 1}R^r)
$$
as a map
$$
K(W,1)\rightarrow K(Center \, \pi _1(R, r) , 2),
$$
in other words as a class in $H^2(K(W,1), Center)$ where $Center$ is the
center of $\pi _1(R, r)$. This is the secondary class we will be interested in
calculating below.

We can describe the secondary class a bit more concretely in the following way.
Suppose $\alpha , \beta$ are elements of $W$ thought of as paths in $K(W,1)$.
Choose trivializations of the fibration
$$
R\rightarrow A\rightarrow K(W,1)
$$
above the paths $\alpha$ and $\beta$. These trivializations lead in particular
to liftings of our paths starting with the basepoint $r$, and ending at points
we denote $\alpha ^{\ast}(r)$ and $\beta ^{\ast}(r)$. By the hypothesis of
vanishing of primary classes, $\alpha ^{\ast}(r)$ and
$\beta ^{\ast}(r)$ are homotopic to $r$. Thus we can choose homotopies which we
denote $h_{\alpha}$ and $h_{\beta}$ (these are paths in $R$ joining
$\alpha ^{\ast}(r)$ resp. $\beta ^{\ast}(r)$ to $r$). Applying the
trivializations chosen above to these paths we obtain paths in $R$
$$
\alpha ^{\ast}(h_{\beta}): \alpha^{\ast} \beta^{\ast} (r) \rightarrow
\alpha ^{\ast} (r),
$$
$$
\beta ^{\ast}(h_{\alpha}): \beta^{\ast} \alpha^{\ast} (r) \rightarrow
\beta ^{\ast} (r).
$$
Finally, the commutativity of $W$ means that there is a torus obtained by
attaching a $2$-cell along the commutator of $\alpha$ and $\beta$. Lifting this
$2$-cell (more precisely, trivializing the family above this $2$-cell)
provides a
path between $\alpha^{\ast} \beta^{\ast} (r)$ and $\alpha^{\ast} \beta^{\ast}
(r)$. Combining these all together we get a $5$-sided loop based at $r$; this is
the secondary class evaluated on $\alpha \wedge \beta$. The loop is in the
center of $\pi _1(R,r)$ because of the fact that we chose trivializations of
the family over our paths (or $2$-cell)
rather than just liftings starting at the basepoint.

A slightly different and more geometric way of looking at this is to look at
the torus in $K(W,1)$ given by the commutator $2$-cell for $\alpha$ and
$\beta$. Attaching two $2$-cells to the torus, one along $\alpha$ and the other
along $\beta$, gives a $2$-sphere. The vanishing of the primary classes
means that the fibration $A\rightarrow K(W,1)$ can be extended across these new
$2$-cells, so we get a family over $S^2$. Trivializing over the northern and
southern hemispheres, the family is determined by a morphism from the equator
to $Aut(R)$: this element of $\pi _1(Aut(R))=Center$ is the image of $\alpha
\wedge \beta$ under the secondary class.

\numero{The secondary Kodaira-Spencer map}

We now come to the situation which gives a ``secondary Kodaira-Spencer map''.
Suppose $f:X_S\rightarrow S$ is a smooth projective morphism. Fix a basepoint
$s\in S$ and denote by $X$ the fiber of $X_S$ over $s$.
Define
$$
W:= T(S)_s
$$
which is a vector space, thus unipotent abelian group scheme, considered as a
sheaf of groups on $Sch /\cc $.

We have an exact
sequence (of unipotent abelian group schemes i.e. vector bundles over $X$)
$$
0\rightarrow  TX \rightarrow T(X_S)|_X \rightarrow f^{\ast}(T(S)_s)\rightarrow
0.
$$
This is an action of $W=T(S)_s$ on the stack $K(TX/X, 1)= X_{UDol}$. In
particular, for any $n$-stack $T$ we obtain an action of $W$ on
$Hom (X_{UDol}, T)$.

The primary invariant in this situation is an action of $W$ on the sheaf of
sets \linebreak $\pi _0Hom (X_{UDol}, T)$.  For example,
if $T=K(\Oo , n)$ then by Proposition \ref{calcDol}
$$
\pi _0Hom (X_{UDol}, T)= H^n_{Dol}(X)= \bigoplus _{p+q=n}H^q(X, \Omega ^p_X)
$$
and we obtain an action of $W=T(S)_s$ on $H^n_{Dol}(X)$. This action is of
course just the usual Kodaira-Spencer map which decomposes into components
$$
T(S)_s \rightarrow Hom \left(  H^q (X, \Omega ^p_X), H^{q+1}(X, \Omega
^{p-1}_X) \right) .
$$

To obtain secondary invariants, we proceed as described above, using vanishing
of the primary invariants if we want to
(but bearing in mind that the secondary invariants will then only be defined
when the primary invariants vanish).  For example, suppose that $T$ is an
$n$-stack and suppose $\eta \in Hom (X_{UDol}, T)$ such that the point $[\eta ]
\in \pi _0Hom (X_{UDol}, T)$ is fixed by the action of $W$. Then $W$ acts on the
connected $n$-stack $Hom ^{\eta}(X_{UDol}, T)$ which is the connected component
containing $\eta$.  We obtain a morphism
$$
K(W,1)\rightarrow B\, Aut (Hom ^{\eta}(X_{UDol}, T)),
$$
and as remarked above, cohomology classes on the right can be pulled back to
give classes on $K(W,1)$. At this point we refer to the  Breen
calculations \cite{Breen2}.  In Appendix II we prove a relative version in
characteristic zero which was stated in \cite{kobe}; the reader may refer there
for the general statement.
In our case, as $W$ is represented by a
finite-dimensional vector space (in particular, $W\cong \Oo ^a$) we have
$$
H^i(K(W,1), \Oo ) = \bigwedge ^iW^{\ast}
$$
where $W^{\ast} = Hom (W, \Oo )$ is the sheaf represented by the dual vector
space.

To get down to the concrete example we would like to consider, we boil things
down a bit farther, following the discussion at the end of the previous section.
Namely instead of looking at the full  $Hom ^{\eta}(X_{UDol}, T)$ we truncate it
down to a $1$-stack (which is connected, too) by looking at
$$
\tau _{\leq 1}Hom ^{\eta}(X_{UDol}, T).
$$
In the very presentable case, this automatically has a base point over $Spec
(\cc )$ (but also we have chosen a basepoint $\eta$) so it is equivalent to
something of the form $K(G,1)$ with
$$
G= \pi _1(Hom (X_{UDol}, T), \eta ).
$$
Again,
in the very presentable case (i.e. if $T$ is very presentable which implies the
same for $Hom ^{\eta}(X_{UDol}, T)$) then $G$ is an affine algebraic group.
In our
example below $G$ will itself be a vector space. The first invariant is
$$
W\rightarrow Out \left( \pi _1(Hom (X_{UDol}, T), \eta )\right) .
$$
Suppose that this
invariant vanishes. Then we get  a map
$$
K(W,1)\rightarrow K(Center \, \pi _1(Hom (X_{UDol}, T), \eta ) , 2),
$$
which may be interpreted as a class in $H^2(K(W,1), Center)$ where $Center$ is
the   center of $\pi _1(Hom (X_{UDol}, T), \eta )$. If $T$ is very presentable
then  this center will be an affine abelian group scheme. This class in $H^2$ is
the class we are interested in calculating below.

The secondary class can be described concretely by choosing homotopies
trivializing the primary classes and combining them together using the
commutativity homotopy for the action of $W$, as at the end of the previous
section.

\numero{The complexified $2$-sphere}

We discuss in more detail the example of a $3$-stack $T$ for which we make our
calculation. Recall from ``standard topology''  how to describe the
rational homotopy type of $S^2$. The only nontrivial stages in the Postnikov
tower are $K(\qq , 2)$ and $K(\qq , 3)$. Thus the rational homotopy type is
described by the fibration sequence
$$
K(\qq , 3) \rightarrow S^2 \otimes \qq \rightarrow K(\qq , 2).
$$
In turn, the classifying space for fibrations with fiber $K(\qq , 3)$ is the
base for the universal fibration
$$
K(\qq , 3)\rightarrow \ast \rightarrow K(\qq , 4).
$$
Thus the above fibration is determined by a morphism
$$
K(\qq , 2)\rightarrow K(\qq , 4),
$$
in other words a class  $\sigma  \in H^4(K(\qq , 2), \qq )$.
The classical calculations give
$$
H^{2m}(K(\qq , 2), \qq ) = Sym ^m(\qq ).
$$
In particular
$$
H^4(K(\qq , 2), \qq )\cong Sym ^2(\qq ) = \qq .
$$
Up to change of basis element for $\pi _3$, there are only two possibilities:
either $\sigma = 0$ or $\sigma \neq 0$. The case $\sigma = 0$ corresponds to
the direct product $K(\qq , 2)\times K(\qq , 3)$ but in this case $H^3$
would be nonzero, whereas $H^3(S^2\otimes \qq , \qq )=0$. Therefore we must
have $\sigma \neq 0$, in other words $\sigma$ is the cup product
$\eta \cup \eta$ of the canonical class $\eta \in H^2(K(\qq , 2), \qq )$
with itself.

In view of this discussion, we can do the same with very presentable stacks. We
will construct $T= S^2\otimes \cc$ with two stages $K(\Oo , 2)$ and $K(\Oo ,
3)$ in the Postnikov tower. We construct it as the pullback of the universal
fibration
$$
K(\Oo , 3)\rightarrow \ast \rightarrow K(\Oo , 4)
$$
by  a morphism
$$
\sigma : K(\Oo , 2) \rightarrow K(\Oo , 4).
$$
The Breen calculations say that the cohomology of $K(\Oo , n)$
with coefficients in $\Oo$ has the same answer as the cohomology of $K(\qq , n)$
with coefficients in $\qq$. In other words,
$$
H^{2m} (K(\Oo , 2), \Oo ) = Sym ^m_{\Oo} (\Oo ) \cong \Oo
$$
and
$$
H^4(K(\Oo , 2), \Oo )= \Oo .
$$
Let $\sigma$ be the generator of $\Oo (Spec \, \cc ) = \cc$. This gives a map
$$
\sigma : K(\Oo , 2)\rightarrow K(\Oo , 4).
$$
In terms of cohomology operations, a morphism $\Ff \rightarrow K(\Oo , 2)$
corresponds to a cohomology class $\eta \in H^2(\Ff , \Oo )$ and the
composition of such a map with $\sigma$ corresponds to the cup-product square
$\eta \cup \eta \in H^4(\Ff , \Oo )$.

Now set
$$
T:= S^2\otimes \cc := K(\Oo , 2) \times _{K(\Oo , 4)} \ast ,
$$
(this choice of $T$ shall be in vigor for the rest of the paper unless
explicitly mentionned otherwise). For any $n$-stack $\Ff$, a morphism
$$
\Ff \rightarrow T
$$
corresponds to a pair $(\eta , \varphi )$ where $\eta : \Ff \rightarrow K(\Oo ,
2)$ and $\varphi$ is a homotopy between
$$
\sigma \circ \eta = \eta \cup \eta : \Ff \rightarrow K(\Oo , 4)
$$
and the constant map at the basepoint
$$
\underline{o}: \Ff \rightarrow K(\Oo , 4).
$$
Alternatively, $\varphi$ may be thought of as a section of the fibration
$$
K(\Oo , 3) \rightarrow \Ff \times _{K(\Oo , 4)} \ast \rightarrow \Ff .
$$
A first remark is that for a given class $\eta \in H^2(\Ff , \Oo )$ there
exists a lifting to a map $\Ff \rightarrow T$ if and only if $\eta \cup \eta =
0$ in $H^4(\Ff , \Oo )$.  If this cup-product is zero so that there exists one
lifting, then the fiber of the map
$$
Hom (\Ff , T) \rightarrow Hom (\Ff , K(\Oo , 2))
$$
over the point $\eta$ is equivalent to $Hom (\Ff , K(\Oo , 3))$
This is because if there exists one  lifting then we can choose a lifting to
trivialize the fibration
$$
\Ff \times _{K(\Oo , 4)} \ast \rightarrow \Ff
$$
i.e. to make this fibration equivalent to
$$
\Ff \times K(\Oo , 3) \rightarrow \Ff
$$
and then the stack of other liftings is just $Hom (\Ff , K(\Oo , 3))$.

{\bf Example:} We give an example of a topological space $Y$ which admits no
nonconstant maps to the actual $2$-sphere $S^2$ but which admits maps to the
complexified $2$-sphere $S^2\otimes \cc$.  Construct $Y$ by taking a wedge of
two $2$-spheres and adding on a $4$-cell via an attaching map
$$
f: S^3\rightarrow S^2 \vee S^2.
$$
Note that
$$
\pi _3(S^2 \vee S^2)= Sym ^2 \pi _2(S^2\vee S^2) = \zz ^3,
$$
and the class of the attaching map $f$ determines the cup product. We can think
of the class of $f$ as a symmetric $2\times 2$ matrix, which can be chosen
arbitrarily. Choose the matrix to be diagonal with $(r,s)$ on the diagonal.
The two obvious classes $e,f\in H^2(S^2\vee S^2,\qq )$ persist as classes in
$H^2(Y,\qq )$. Note that $H^4(Y,\qq ) \cong \qq $ and we have the formulae
$$
e\cup e = r, \;\; f\cup f = s,\;\; e\cup f = 0.
$$
Now a map $Y\rightarrow S^2$ (or even to $S^2\otimes \qq$) corresponds to a
class
$$
\eta = ae +bf \in H^2(Y,\qq )
$$
with $a,b\in \qq $, such that $\eta \cup \eta =0$. The lifting to a map into
$S^2$ is unique because in our case $H^3(Y,\qq )=0$. However, if $r$ and $s$
are chosen to be relatively prime and having no square prime factors,
then the equation
$$
\eta \cup \eta = a^2r + b^2s = 0
$$
doesn't have any nonzero solutions with $a,b\in \qq$ (the same holds if $r$ and
$s$ have the same sign).  Thus there are no nontrivial maps $Y\rightarrow S^2$.
On the other hand, the above equation always has nonzero complex solutions,
so there is always a nontrivial map $Y\rightarrow S^2\otimes \cc $.

\numero{Our example}
We now come down to the example which we would like to calculate. Fix $T=
S^2\otimes \cc$ as defined above. Let $Z$ be a smooth projective surface with
$H^1(Z, \Oo ) = 0$ (hence $H^1(Z, \cc )= H^3(Z, \cc )=0$).  Let $P\in Z$ and
let $X$ be the blow-up of $Z$ at $P$. Let $W\subset H^1(X, TX)$ be the rank two
subspace of deformations of $X$ corresponding to moving the point $P$ which is
blown up (thus canonically $W\cong T(Z)_P$). Then $W$ acts on $X_{UDol}$, via
the exact sequence
$$
0\rightarrow TX\rightarrow T(Tot)|_X \rightarrow W\otimes _{\cc}\Oo _X
\rightarrow 0
$$
where $Tot$ refers to the total space of the family.
Let
$$
R:= Hom (X_{UDol}, T).
$$
Then $W$ acts on $R$ and we will look at secondary classes for this action.

Because of our hypothesis $H^3_{Dol}(X)=0$, the map $T\rightarrow K(\Oo , 2)$
induces an injection
$$
\pi _0(R) \hookrightarrow H^2_{Dol}(X)
$$
with the image the sheaf represented by the subscheme defined by the equation
$\eta \cup \eta = 0$. The variation of Hodge structure of $H^2(X, \cc )$
parametrized by $P\in Z$ as we move the point which is blown up, is trivial.
Thus the Kodaira-Spencer class is trivial, in other words $W$ acts trivially on
$H^2_{Dol}(X)$ and hence it acts trivially on $\pi _0(R)$.  If we fix
$\rho = (\eta , \varphi )\in R(Spec \, \cc )$ then we obtain a secondary class
$$
\kappa \in H^2(K(W, 1), \pi _1(R, \rho )).
$$
The first task is to calculate $\pi _1(R, \rho )$ (which will be a group scheme
over $Spec (\cc )$ since $\rho$ is defined over $Spec (\cc )$). Recall
that we have the fibration
$$
R\rightarrow Hom (X_{UDol}, K(\Oo , 2))
$$
whose fiber over a point in the target is either empty or else equivalent to
\newline
$Hom (X_{UDol}, K(\Oo , 3))$.  Thus the long exact sequence for this fibration
gives
$$
\cc = H^0_{Dol}(X) \rightarrow H^2_{Dol}(X)\rightarrow \pi _1(R, \rho )
\rightarrow H^1_{Dol}(X) = 0.
$$
the first term being $\pi _2Hom (X_{UDol}, K(\Oo , 2))$, the
second term being
\newline
$\pi _1Hom (X_{UDol}, K(\Oo , 3))$, and so on.

{\bf Claim:} The connecting morphism $H^0_{Dol}(X) \rightarrow H^2_{Dol}(X)$
in the above exact sequence is multiplication by a nonzero multiple of $\eta$.

{\em Proof:} For any $n$-stack $Y$ we can define $Hom (Y, T)$ and look at the
long exact sequence for the fibration
$$
Hom(Y,T)\rightarrow Hom(Y,K(\Oo , 2)).
$$
The connecting morphism will be functorial in $Y$. Apply this to $Y=T$ itself;
then for any other $Y$ (such as $X_{UDol}$ considered in the claim) the
connecting morphism for the long exact sequence at a basepoint $Y\rightarrow T$
is obtained by pulling back the connecting morphism for $Hom (T,T)$ at the
identity map.
In other words, looking at $Y=T$ gives a universal version of the connecting
morphism. On the other hand, note that
$$
\pi _1 Hom (T, K(\Oo , 3))=H^2(T, \Oo )=\Oo
$$
and
$$
\pi _2Hom (T,K(\Oo , 2)) = H^0(T,\Oo )=\Oo .
$$
Thus the universal connecting morphism is a scalar constant $C$, and for any
$\rho : Y\rightarrow T$ the connecting morphism for $Hom (Y,T)$ based at $\rho$
fits into the diagram
$$
\begin{array}{ccc}
H^0(T,\Oo ) & \rightarrow & H^2(T, \Oo )\\
\downarrow && \downarrow \\
H^0(Y,\Oo ) & \rightarrow & H^2(Y, \Oo ).
\end{array}
$$
The vertical maps are those induced by $\rho$. It follows that the connecting
map for $Hom (Y,T)$ based at $\rho$ is the same constant $C$ multiplied by the
class $\eta$ which is the image of $\rho$ in $H^2(Y,\Oo )$. For the claim, apply
this to $Y=X_{UDol}$. To finish proving the claim we just have to show that $C$
is nonzero. But if $C$ were zero then the generator for $\pi _2(Hom (T,K(\Oo ,
2)), 1_T)$ would lift to  an element of $\pi _2(Hom(T,T), 1_T)$. This would give
a map  $$
S^2 \times T \rightarrow T
$$
(where $S^2$ denotes the constant presheaf with values $S^2$),
which is nontrivial on $S^2\times \{ 0\} $ and $\{ 0\} \times T$. This map
would correspond to an element
$\mu \in H^2(S^2\times T, \Oo )$ with $\mu \cup \mu = 0$. But we have
$$
H^2(S^2\times T, \Oo )= H^2(S^2, \Oo )\oplus H^2(T, \Oo ) = \Oo \oplus \Oo
$$
and the cup product of the two components is nontrivial (by K\"unneth).
Therefore it is impossible to have a class $\mu$ which is nontrivial in both
components but with $\mu \cup \mu = 0$. This contradiction implies that $C\neq
0$, giving the claim.
\eop

With the claim we obtain
$$
\pi _1(R,\rho ) \cong H^2_{Dol}(X)/(\eta ).
$$
Thus our characteristic class $\kappa $ becomes a map
$$
\kappa : \bigwedge ^2(W) \rightarrow H^2_{Dol}(X)/(\eta ).
$$

Next, we choose $\eta$ (which fixes the choice of $\rho$ up to homotopy).
Let $E$ be the exceptional divisor on $X$ and let $H$ denote the pullback of an
ample divisor on $Z$ not meeting the point $P$. Let $[E]$ and $[H]$ denote
their Chern classes in $H^1(X, \Omega ^1_X)\subset H^2_{Dol}(X)$.
We set
$$
\eta = m[E] + n[H] \in H^2_{Dol}(X).
$$
We have to choose $n$ and $m$ so that $\eta \cup \eta = 0$. Note that
$H^4_{Dol}(X)\cong \cc $ with natural morphism given by the residue map,
normalized so that the cohomology class of a point is equal to $1$. Via this
isomorphism, $[E] \cup [E] = E.E = -1$ and $[H] \cup [H] = H.H\in \zz $.
Note that $[E]\cup [H]=0$ since the two divisors don't intersect.
Thus
$$
\eta \cup \eta = n^2H.H - m^2.
$$
We choose $m,n$ so that this is equal to $0$.

{\bf Remark:} If $H.H$ is not  the square of an integer, then the pair
$(m,n)$ can not be chosen in $\qq ^2$, and in particular our map will not exist
as a topological map $X^{\rm top} \rightarrow S^2$.  However, our map will
exist as a map from the constant
presheaf with values $X^{\rm top}$, to $S^2\otimes \cc $. (This is a heuristic
remark since in the present paper we don't treat the question of the
relationship between Betti cohomology and Dolbeault cohomology).

Here is the result of the calculation which will be done below.

\begin{theorem}
\mylabel{calculation}
With the above choices of $T$, $X$,
$R:= Hom (X_{UDol}, T)$, $\eta = m[E] + n[H]$ ($\eta \neq 0$), and
hence $\rho \in R(Spec \, \cc )$, the secondary Kodaira-Spencer class
$$
\kappa : \bigwedge ^2 W \rightarrow \pi _1(R, \rho )=
H^2_{Dol}(X)/(\eta )
$$
lands in $H^2(X, \Oo ) \subset H^2_{Dol}(X)/(\eta )$ and the map
$$
\bigwedge ^2 W \rightarrow H^2(X,\Oo )
$$
is dual (using Serre duality,
the isomorphism $H^0(X, \Omega ^2_X)\cong H^0(Z, \Omega ^2_Z)$, and
the isomorphism $W^{\ast} \cong T^{\ast}(Z)_P$)
to $m^2ev_P$ where
$$
ev_P:H^0(Z, \Omega ^2_Z) \rightarrow \bigwedge ^2T^{\ast}(Z)_P
$$
is the evaluation morphism.
In particular, $\kappa \neq 0$ if $h^{2,0}(X)=h^{2,0}(Z) > 0$ and $m\neq 0$.
\end{theorem}

Before getting on with the proof in \S 8 below, we make a few general remarks
about this result. A similar thing can be stated for the ``Dolbeault
homotopy type of $X$''. One way of defining this (which wouldn't be the
historical way, though) is as the $1$-connected very presentable $n$-stack
$\Sigma$
representing the very presentable  shape of $X_{Dol}$ (cf Theorem
\ref{representable1} in Appendix II below).  In this point of view,
we get an action of $W$ on the Dolbeault homotopy type.  The theorem says
that this action of $W$ is nontrivial.
Note however that the action of $W$ on the homotopy group sheaves
(which are the homotopy groups of $X$ tensored with $\cc$) will be
trivial. It is certainly possible to define the action of $W$, and to
make the same calculation as below to show that the action is nontrivial, using
the algebra of forms  $$
A^{\cdot}
_{Dol}(X)= (\bigoplus _{p,q} A^{p,q}(X), \delbar ).
$$
In fact, my first heuristic version of the calculation was done using forms.
However, the technical details relating a differential-forms version of
nonabelian cohomology, with the version presented here, seem for the moment
somewhat difficult, so we restrict in the present paper to an algebraic
version of the calculation.

The secondary class $\kappa$ is a natural map, so it doesn't really have any
choice other than to be a multiple of the dual of the evaluation map $ev_P$. The
only question is whether this multiple is nonzero or not.  Here is a heuristic
global argument to see why, in principle, the constant should be nonzero.
Let $\Xx \rightarrow Z$ be the total space of the family of blow-ups of points
moving in $Z$. It is obtained by blowing up the diagonal in $Z\times Z$. The
secondary Kodaira-Spencer class we calculate here is (or should be, at least)
the $(2,0)\times (0,2)$ Hodge component
(i.e. the component of type $(2,0)$ on the base and $(0,2)$ on the fiber) of the
following globally defined invariant. Fix $\eta \in H^2(X,\cc )$, which is
invariant under the monodromy since the monodromy is trivial (the base $Z$ being
simply connected); the degeneration of the Leray spectral sequence says that
$\eta$ comes from restriction of a global class $\tilde{\eta} \in H^2(\Xx , \cc
)$; the cup product $\tilde{\eta} \cup \tilde{\eta}$ restricts to zero on the
fibers,  so it lies in the next step for the Leray filtration, which in our case
is $H^2(Z, H^2(X, \cc ))$. This cup product is therefore a globally defined
class. It is the obstruction to extending $\rho : X^{\rm top}\rightarrow T$ to a
map $\Xx ^{\rm top}\rightarrow T$. The $(2,0)\times (0,2)$ component of this
class, which is a holomorphic $2$-form on $Z$ with coefficients in
$H^2(X, \Oo )$, {\em should} give $\kappa$ when evaluated at the point $P\in Z$.
I don't currently have a proof of this, though.

In our example, it is relatively easy to see by looking at the
cohomology class of the diagonal that the
$(2,0)\times (0,2)$ component of the global class is
nonzero, and in fact it is the identity matrix (via Serre duality).
Thus if one could prove the above statement that the global class gives $c$
under evaluation at $P\in Z$, then this would prove that $\kappa\neq 0$ when
$H^2(X, \Oo )\neq 0$.

With the previous paragraphs as heuristic argument, the result that $\kappa\neq
0$ doesn't look all that surprising. Still, it means that the ``variation of
nonabelian Hodge structure''
\footnote{
This terminology is put in quotes because the current discussion of Dolbeault
cohomology is only a first step towards defining what a
``variation of
nonabelian Hodge structure'' is.}
on the family of $Hom (X^{\rm top}, T)$,
when $X$ is a variable fiber in the family $\Xx \rightarrow Z$, is nontrivial,
and even nontrivial for infinitesimal reasons. The base of this ``variation'' is
$Z$ which is simply connected.  In particular, the variations of mixed Hodge
structure on the homotopy groups (or anything else you could think of) are
trivial. From a topological point of view, it is never a surprise to find a
family where the homotopy groups are constant but the family nontrivial. On the
other hand, this goes against the commonly held intuition for
projective algebraic varieties that ``formality means that everything is
determined by the cohomology ring'': in the example $\Xx \rightarrow Z$, the
locally constant family---parametrized by $Z^{\rm top}$---of cohomology rings of
the fibers $X$ is trivial. What remains true of course is that the topology of
the family is determined by the cohomology ring of the total space $\Xx$. Our
secondary Kodaira-Spencer invariant is  a local invariant which
contributes to nontriviality of the global cohomology ring of the total space of
the family.

Our class detects the motion of a point $P\in Z$ exactly when $H^0(Z,\Omega
^2_Z)\neq 0$.  This seems to fit in with the standard intuition that
$H^0(Z,\Omega ^2_Z)\neq 0$ causes the class group of zero cycles to be big
(Mumford's and Clemens' results, Bloch conjecture etc. cf
Voisin \cite{Voisin}).  I don't see a precise connection, though.

\numero{The calculation}

We keep the above notations $T$, $P\in Z$, $X$, $E$, $H$, $\eta$, $\rho$.
We establish some more: let $N$ be an affine neighborhood of $P$ in $Z$
such that $TZ$ is trivialized over $N$. Assume that $N$ doesn't meet the
divisor image of $H$ in $Z$. Let $\alpha$ and $\beta$ denote basis sections in
$TZ(N)$. Let $B$ be the inverse image of $N$ in $X$, and let $C=X-E \cong Z-\{
P\}$. Then put $A:= B\cap C$. Note that $\{ B,C\}$ is an open covering of $X$.
We can write $$
X= B\cup ^AC
$$
(this is true as a pushout of sheaves of sets on $Sch/\cc $).  Similarly we have
$$
X_{UDol}= B_{UDol} \cup ^{A_{UDol}}C_{UDol}.
$$
The basis vectors $\alpha , \beta$ give cocycles for elements in $H^1(X, TX)$
(actually in the \v{C}ech cohomology relative to our covering) and it is easy
to see that these cocycles project to basis elements of our
$2$-dimensional space $W$.
We denote the basis vectors of $W$ also by $\alpha$ and $\beta$.

We can describe the action of $W$ on $X_{UDol}$ concretely in the following
way. Set $K:= K(W, 1)$ with basepoint denoted $0\in K$. Then the trivialization
of $TZ|_N$ gives an equivalence
$$
A_{UDol} \cong A\times K.
$$
There is a group structure on $K$, that is a morphism $K\times K\rightarrow K$
corresponding to the addition on $W$, and this gives
$$
A\times K\times K\rightarrow A\times K,
$$
which we can rewrite as
$$
\mu : A_{UDol}\times K \rightarrow A_{UDol}.
$$
Putting the identity $1_K$ in the second variable we get a map
$$
\Phi: A_{UDol} \times K \rightarrow A_{UDol} \times K
$$
which is an equivalence. Note also that
$\Phi |_{A_{UDol} \times 0}$ is the identity of $A_{UDol}$.

Let $i,j$ be the
inclusions from $A_{UDol}$ to $B_{UDol}$ and $C_{UDol}$ respectively. Then use
the inclusions
$$
(i\times 1_K)\circ \Phi : A_{UDol}\times K \rightarrow B_{UDol}\times K
$$
and
$$
j\times 1_K: A_{UDol}\times K \rightarrow C_{UDol}\times K
$$
to construct the pushout
$$
P:= B_{UDol}\times K\cup ^{A_{UDol}\times K}C_{UDol}\times K.
$$
This comes equipped with a morphism
$$
P\rightarrow K
$$
and the fiber over the basepoint $0\in K$ is just
$$
B_{UDol}\cup ^{A_{UDol}}C_{UDol} = X_{UDol}.
$$
Therefore according to our definition, $P\rightarrow K$ is an action
of $W$ on $X_{UDol}$.

The corresponding action of $K$ on $R:= Hom (X_{UDol}, T)$ is by definition
$$
Hom (P/K, T)\rightarrow K.
$$
Using the Mayer-Vietoris principle we get
$$
Hom (P/K, T) =
$$
$$
Hom (B_{UDol}\times K/K,T) \times _{Hom (A_{UDol}\times K/K,T)}
Hom (C_{UDol}\times K/K,T).
$$
However, note that
$$
Hom (B_{UDol}\times K/K,T) = Hom (B_{UDol}, T) \times K
$$
and similarly for the other factors. The morphism
$$
Hom (C_{UDol}, T) \times K\rightarrow
Hom (A_{UDol}, T) \times K
$$
induced by $j\times 1_K$ is just the product of the morphism induced by $j$,
with $1_K$.  (On the other hand, the same is not true of the first morphism in
the fiber product, as  it is induced by $(i\times 1_K)\circ \Phi$.) We can now
write
$$
Hom (P/K, T) =
$$
$$
(Hom (B_{UDol},T)\times K) \times _{Hom (A_{UDol},T)}
Hom (C_{UDol},T).
$$
The first morphism in the fiber product is the composition of the
product-compatible morphism
$$
(j^{\ast} \times 1_K): Hom (B_{UDol},T)\times K
\rightarrow Hom (A_{UDol},T)\times K,
$$
with the morphism
$$
\Psi : Hom (A_{UDol},T)\times K\rightarrow Hom (A_{UDol}, T).
$$
This map is equivalent (by the definition of internal $Hom$) to
$$
Hom (A_{UDol},T)\rightarrow Hom (K, Hom (A_{UDol},T)),
$$
which in turn is equivalent to a map
$$
Hom (A_{UDol},T)\rightarrow Hom (A_{UDol}\times K, T),
$$
this latter being induced by our action $\mu : A_{UDol}\times K \rightarrow
A_{UDol}$.

The first step is to notice that the map $\rho |_{A_{UDol}}$ from $A_{UDol}$ to
$T$ factors through a map
$$
h: A_{UDol} \rightarrow K(\Oo , 3)\rightarrow T.
$$
This factorization is given by the fact that $\eta$ is a class in
$H^2((X_{Dol}, A_{Dol}), \Oo )$, in other words we are given a trivialization
of $\eta$ over $A_{Dol}$.

Recall that we write $\rho = (\eta , \varphi )$ where $\eta :
X_{UDol}\rightarrow K(\Oo , 2)$ and $\varphi$ is a section of the pullback
bundle
$$
L_{\eta}:= X_{UDol} \times _{K(\Oo , 4)} \ast
$$
which is a bundle with fiber $K(\Oo , 3)$ over $X_{UDol}$. The section
$\varphi$ determines a trivialization
$$
L_{\eta} \cong X_{UDol}\times K(\Oo , 3)
$$
such that $\varphi$ corresponds to the $0$-section. This trivialization is
uniquely determined by the condition that it be compatible with the structure
of ``principal bundle'' under the ``group'' $K(\Oo , 3)$.

We adopt the following strategy for calculating the secondary class. We
will look
at new $n$-stacks $P^i\rightarrow K$ with $K$-maps $P^i \rightarrow P$, and
points $\rho ^i\in P^i_0$ (where $P^i_0$ means the fiber of $P^i$ over
$0\in K$), such that $\rho ^i$ maps to $\rho$. We have to arrange so that
the primary class is trivial, in other words that the class of $\rho ^i$ in $\pi
_0(P^i_0)$ should be invariant under the action of $W=\pi _1(K)$.
We also have to insure that the other primary class, the action of $W$ on $\pi
_1(P^i_0, \rho ^i)$ by outer automorphisms, should be trivial. In this case,
we obtain a secondary class for $\rho ^i$, which is an element of
$$
H^2(K, \pi _1(P^i_0, \rho ^i))
$$
and this secondary class maps to our class for $P$.

First of all, let $\ast \rightarrow Hom (B_{UDol}, T)$ be the morphism
corresponding to the point $\rho |_{B_{UDol}}$.  We get
$$
K= \ast \times K \rightarrow Hom (B_{UDol}, T)\times K.
$$
Thus we obtain a morphism
$$
P^1:= K\times _{Hom (A_{UDol}, T)}Hom (C_{UDol}, T) \rightarrow
$$
$$
(Hom (B_{UDol}, T)\times K)
\times _{Hom (A_{UDol}, T)}Hom (C_{UDol}, T)
=P.
$$
This morphism is compatible with the projections to $K$.

The first morphism
$$
u: K\rightarrow Hom (A_{UDol}, T)
$$
in the fiber product is obtained by the composition
$$
K \rightarrow Hom (B_{UDol}, T)\times K\stackrel{(i\times  1_K)\circ
\Phi}{\rightarrow}
$$
$$
Hom (A_{UDol}, T)\times K \stackrel{p_1}{\rightarrow}
Hom (A_{UDol}, T).
$$
Thus $u$ corresponds to the map $K\times A_{UDol} \rightarrow T$ obtained by
composing
$$
K\times A_{UDol} \stackrel{\mu }{\rightarrow}
A_{UDol} \stackrel{\rho}{\rightarrow} T.
$$
The map $u$
factors through a morphism
$$
K\rightarrow \Gamma (A_{UDol}, L_{\eta})\rightarrow Hom (A_{UDol}, T)
$$
(technically speaking what should enter into the above
notation is $L_{\eta}|_{A_{UDol}}$ but for brevity we omit the restriction
since it is implicitly determined by the notation $\Gamma (A_{UDol}, -)$.)
The factorization comes about because  the composition
$$
A_{UDol} \stackrel{\rho}{\rightarrow} T \rightarrow K(\Oo , 2),
$$
is given as the constant map at the basepoint. Thus pulling back by
$$
A_{UDol}\times K \rightarrow A_{UDol} \rightarrow K(\Oo , 2)
$$
is again constant at the basepoint so the map
$$
K\rightarrow Hom( A_{UDol}, T)
$$
factors through a map
$$
\tilde{u}:K\rightarrow Hom( A_{UDol}, K(\Oo , 3)).
$$
Following above, $\tilde{u}$ corresponds to the composition
$$
K\times A_{UDol}\stackrel{\mu}{\rightarrow} A_{UDol} \stackrel{h}{\rightarrow }
K(\Oo , 3),
$$
which we can write as $\tilde{u} = \mu ^{\ast}(h)$.

Now we obtain a morphism (over $K$)
$$
P^2:= K\times _{\Gamma (A_{UDol}, L_{\eta})}\Gamma (C_{UDol}, L_{\eta})
$$
$$
\rightarrow
K\times _{Hom (A_{UDol}, T)}Hom (C_{UDol}, T) = P^1.
$$
We have used the section $\varphi$ of $L_{\eta}$ to obtain a trivialization
$$
L_{\eta} \cong X_{UDol} \times K(\Oo , 3).
$$
Via this equivalence
the section corresponding to $\rho$ (that is, the
section $\varphi$) corresponds to the zero-section.
Using the trivialization given by $\varphi$ we can write
$$
\Gamma (A_{UDol}, L_{\eta})\cong Hom (A_{UDol}, K(\Oo , 3))
$$
({\em Caution:} this is not the same trivialization as given by saying that
$A_{UDol}$ maps to the basepoint of $K(\Oo , 2)$ so the pullback fibration
$L_{\eta}$ is trivial over $A_{UDol}$; these two trivializations differ by
translation by $h$, a point which will come up below);
and
$$
\Gamma (C_{UDol}, L_{\eta})\cong Hom (C_{UDol}, K(\Oo , 3)).
$$
These are compatible with the restriction $j^{\ast}$ so we can write
$$
P^2 = K\times _{Hom (A_{UDol}, K(\Oo , 3))} Hom (C_{UDol}, K(\Oo , 3)),
$$
with the second morphism in the fiber product being the restriction $j^{\ast}$
acting on maps to $K(\Oo , 3)$.

We have to re-calculate the first morphism in the fiber product
$$
a:K\rightarrow Hom (A_{UDol}, K(\Oo , 3)).
$$
It is no longer equal to $\tilde{u}=\mu ^{\ast}(h)$, because when we set
the section $\varphi$ of $L_{\eta}$ equal to the zero-section to get an
equivalence between $L_{\eta}$ and $K(\Oo , 3)$, this made a translation on
$L_{\eta} |_{A_{UDol}}$---which was already trivial due to the fact that $\eta
|_{A_{UDol}}=0$---this translation has the effect of setting $h$ equal to the
$0$-section.  This translation gives us the formula
$$
a = \mu
^{\ast}(h)-p_2^{\ast}(h) $$ where $p_2^{\ast}(h)$ is the map $K\rightarrow Hom
(A_{UDol}, K(\Oo , 3))$ corresponding to the composition
$$
K\times A_{UDol} \stackrel{p_2}{\rightarrow}
A_{UDol} \stackrel{h}{\rightarrow} K(\Oo , 3).
$$
Note that $p_2$ denotes the second projection. The minus sign in the equation
for $a$ is subtraction using the ``abelian group'' (i.e. $E_{\infty}$) structure
of $Hom(A_{UDol}, K(\Oo , 3))$ induced by the ``abelian group'' structure of
$K(\Oo , 3)$.

We now turn back to the second morphism in the fiber product,
$$
j^{\ast} : Hom( C_{UDol}, K(\Oo , 3))\rightarrow Hom( A_{UDol}, K(\Oo , 3)).
$$
The infinite loop space structure of $K(\Oo , 3)$ which is inherited
by  $Hom (\Ff , K(\Oo , 3))$ for  $\Ff = C_{UDol}$ and $\Ff = A_{UDol}$.
By Proposition \ref{decomp}, this delooping structure gives  a
decomposition of
$Hom (\Ff , K(\Oo , 3))$ into a product of Eilenberg-MacLane stacks.
The restriction morphism $j^{\ast}$ above is compatible with the
delooping structures, so by Proposition \ref{decomp} (B), it is homotopic
to a map
compatible with the decomposition into a product of Eilenberg-MacLane stacks.

Recall that
$$
\pi _i(Hom (C_{UDol}, K(\Oo , 3)), 0) = H^{3-i}(C_{UDol}, \Oo )=
H^{3-i}_{Dol}(C)
$$
and
$$
\pi _i(Hom (A_{UDol}, K(\Oo , 3)), 0) = H^{3-i}(A_{UDol}, \Oo )=
H^{3-i}_{Dol}(A).
$$

From the splitting given by the delooping structures via Proposition
\ref{decomp} with homotopy of functoriality of part (B) of that proposition
(choose one), we
obtain a homotopy-commutative diagram
$$
\begin{array}{ccc}
K(H^2_{Dol}(C), 1)& \rightarrow & K(H^2_{Dol}(A), 1) \\
\downarrow && \downarrow \\
Hom( C_{UDol}, K(\Oo , 3))&\rightarrow &Hom( A_{UDol}, K(\Oo , 3)).
\end{array}
$$
All maps are infinite loop maps.
Using this diagram, we get the map
$$
P^3:= K\times _{Hom( A_{UDol}, K(\Oo , 3))}K(H^2_{Dol}(C), 1)
\rightarrow
$$
$$
K\times _{Hom (A_{UDol}, K(\Oo , 3))} Hom (C_{UDol}, K(\Oo , 3))=P^2.
$$
Set $\rho ^3:= (0,0)$ in $P^3_0$. It maps to $\rho ^2$ (a point which we didn't
specify because it was always obviously given by $\rho$), because we have
normalized so that $\varphi$ becomes the zero-section. The map $P^3\rightarrow
P^2$ (obviously a map over $K$) takes $\rho ^3$ to $\rho ^2$.

We have
$$
H^2_{Dol}(C)= H^2(C, \Oo _C)\oplus H^1(C, \Omega ^1_C)\oplus H^0(C, \Omega
^2_C),
$$
and similarly
$$
H^2_{Dol}(A)= H^2(A, \Oo _A)\oplus H^1(A, \Omega ^1_A)\oplus H^0(A, \Omega
^2_A).
$$
Recall that $C$ and $A$ are isomorphic to open subsets of $Z$: together $N$ and
$C$ form a covering of $Z$ and $N\cap C= A$. On the other hand,
cohomology of coherent sheaves on $N$ vanishes because $N$ is affine. Therefore
Mayer-Vietoris gives an exact sequence
$$
H^1(C,\Omega ^1_C) \rightarrow H^1(A, \Omega ^1_A) \rightarrow H^2(Z, \Omega
^1_Z).
$$
We are supposing that $H^1_{Dol}(Z)=0$ so by duality $H^3_{Dol}(Z)=0$.
Thus the term on the right is zero, and the morphism of vector spaces
$$
H^1(C,\Omega ^1_C) \rightarrow H^1(A, \Omega ^1_A)
$$
is surjective. Choose a splitting. This gives a morphism
$$
K(H^1(A, \Omega ^1_A), 1)\rightarrow K(H^2_{Dol}(C), 1)
$$
such that the projection into $K(H^2_{Dol}(A), 1)$ is equal to the inclusion
from the Dolbeault direct sum decomposition for $A$. Using our choice of
splitting we obtain a map
$$
P^4:= K\times _{Hom( A_{UDol}, K(\Oo , 3))}K(H^1(A, \Omega ^1_A), 1)
\rightarrow
$$
$$
K\times _{Hom( A_{UDol}, K(\Oo , 3))}K(H^2_{Dol}(C), 1) = P^3.
$$
Again set $\rho ^4 = (0,0)\in P^4_0$, which maps to $\rho ^3$.

Now $P^4$ is defined entirely in terms of $A$.
We decompose things a bit more. Using the infinite loop-space
structure of $Hom( A_{UDol}, K(\Oo , 3))$ and again Proposition \ref{decomp}
we get the decomposition into a product of Eilenberg-MacLane spaces
$$
Hom( A_{UDol}, K(\Oo , 3))\cong
J^0\times J^1 \times J^2 \times J^3
$$
where
$$
J^i = K(H^{3-i}_{Dol}(A), i).
$$
The second map in the fiber product defining $P^4$ is homotopic (choose a
homotopy) to one which is compatible with this decomposition, coming from the
morphism
$$
K(H^1(A, \Omega ^1_A), 1)\rightarrow J^1= K(H^2_{Dol}(A), 1)
$$
(the other components are the maps sending $\ast$ to the basepoints $0\in J^i$,
$i\neq 2$).

Thus we can write
$$
P^4 = Q^0\times Q^1 \times Q^2 \times Q^3
$$
where
$$
Q^i= K\times _{J^i}\ast = K\times _{K(H^{3-i}_{Dol}(A), i)}\ast
$$
for $i\neq 1$ and where
$$
Q^1 = K\times _{K(H^{2}_{Dol}(A), 1)}K(H^1(A, \Omega ^1_A), 1).
$$
This means that the action of $W$ on $P^4_0$ decomposes into a
product of actions
on the $Q^i_0$. Note that the $Q^i_0$ are themselves Eilenberg-MacLane spaces,
$$
Q^i_0 = \Omega J^i = K( H^{3-i}_{Dol}(A), i-1)
$$
for $i\geq 2$, and
$$
Q^1_0 = K(H^{2}_{Dol}(A)/H^1(A, \Omega ^1_A), 0).
$$
The classifying maps for the actions of $W$ on the components $Q^i_0$
factor as
$$
K\rightarrow J^i \rightarrow B\, Aut (Q^i_0)
$$
for $i\geq 2$, and
$$
K\rightarrow K(H^{2}_{Dol}(A)/H^1(A, \Omega ^1_A), 1)\rightarrow
B\, Aut (Q^1_0)
$$
for $i=1$. Thus the classifying map for the action
$K\rightarrow B\, Aut (P^4_0)$
factors through our above map
$$
a: K\rightarrow
K(H^{2}_{Dol}(A)/H^1(A, \Omega ^1_A), 1)\times J^2\times J^3.
$$
We use this to calculate the characteristic classes for the action:
the  primary invariant is the first component, corresponding to a map
$$
W\rightarrow H^{2}_{Dol}(A)/H^1(A, \Omega ^1_A).
$$
The
other primary invariant giving the action of $W$ on $\pi _1$ by outer
automorphisms, is trivial: it is the action on $\Omega J^2$ induced by the
classifying map $K\rightarrow J^2$, and $J^2$ is simply connected.
The secondary invariant which we are interested in is the map
$K\rightarrow J^2$ which corresponds to a class in $H^2(K, H^1_{Dol}(A))$.

A preliminary remark is that the formula $a= \mu ^{\ast}(h)-p_2^{\ast}(h)$
means that the component of $a$ in $J^0$ is equal to $0$. In fact,
$p_2^{\ast}(h)$ is exactly the $J^0$-component of $\mu ^{\ast}(h)$.
On the other hand, the remaining components of $a$ are the same
as those of $\mu
^{\ast}(h)$; these are the K\"unneth components of
$$
h\circ \mu : K\times A_{UDol} \rightarrow K(\Oo , 3).
$$

The first thing to check is that
the primary invariant is trivial for $(P^4, \rho ^4)$.
As we have said above, it is the map
$$
K(W, 1)\rightarrow K(H^{2}_{Dol}(A)/H^1(A, \Omega ^1_A), 1).
$$
To check that it is trivial, it suffices to prove the
\newline
{\bf Claim} ${\bf (\ast )}:$
the map
$$
W\rightarrow \pi _1(Hom( A_{UDol}, K(\Oo , 3)), 0)= H^2_{Dol}(A)
$$
takes $W$ into the component
$H^1(A, \Omega ^1_A)$.  For this we must again get back to the description of
the map
$$
\mu ^{\ast}(h): K\rightarrow Hom (A_{UDol}, K(\Oo , 3))
$$
(note as remarked above that all of the components except the $J^0$ component,
are the same for $\mu ^{\ast}(h)$ or for $a$).
This map, which is equivalent to a
map
$$
A_{UDol}\times K \rightarrow K(\Oo , 3),
$$
is the pullback of
$$
h: A_{UDol}\rightarrow K(\Oo , 3)
$$
by the map
$$
\mu : A_{UDol}\times K \rightarrow A_{UDol}.
$$
Now we have
$$
A_{UDol} = A\times K
$$
so
$$
H^3(A_{UDol}, \Oo ) = H^3(A, \Oo )\oplus H^2(A, H^1(K, \Oo ))
\oplus H^1(A, H^2(K, \Oo ) ) \oplus H^3(K , \Oo ).
$$
The last term $H^3(K , \Oo )$ vanishes because it would be
$\bigwedge ^3W^{\ast}$ but $W$ is $2$-dimensional. Similarly, one can arrange
(by an appropriate choice of $N$) so that $A$ has an open covering by two
affine open sets. Thus $H^3(A, \Oo ) = H^2(A, \Oo )= 0$. The only remaining
term is
$$
H^3(A_{UDol}, \Oo ) = H^1(A, H^2(K, \Oo ))= H^1(A, \bigwedge ^2W^{\ast})
= H^1(A, \Omega ^2_A).
$$
Our map
$$
A_{UDol}\times K \rightarrow K(\Oo , 3)
$$
is obtained by pulling back the above, using the map $K\times K\rightarrow K$.
Pullback for this map is
$$
\bigwedge ^2W^{\ast} \rightarrow
$$
$$
H^2(K,\Oo)\otimes _{\Oo}\Oo \oplus H^1(K,\Oo )\otimes _{\Oo}H^1(K, \Oo ) \oplus
\Oo \otimes _{\Oo}H^2(K,\Oo ).
$$
Each of the factors is nontrivial, with the middle being (up to a multiple
which depends on normalizations for notation in exterior products)
the standard map
$$
\bigwedge ^2 W^{\ast} \rightarrow W^{\ast}\otimes _{\Oo}W^{\ast}.
$$
The morphism induced by pulling back $h$ to $A_{UDol}\times K$ thus decomposes
into K\"unneth components
$$
H^1(A, \Omega ^2_A)\otimes _{\Oo}\Oo \; \; \oplus \; \;
H^1(A, \Omega ^1_A)\otimes _{\Oo}W^{\ast} \; \; \oplus \; \;
H^1(A, \Oo _A)\otimes _{\Oo}\bigwedge ^2W^{\ast}.
$$
Each component is induced by $h\in H^1(A,\Omega ^2_A)$.

The map
$$
\pi _1(K)=W\rightarrow \pi _1(Hom (A_{UDol}, K(\Oo , 3))
= H^2_{Dol}(A)
$$
corresponds to the middle component, which is in fact a map
$$
W\rightarrow H^1(A, \Omega ^1_A).
$$
This was exactly the claim ${\bf (\ast )}$ we needed to prove to show that the
primary invariant for $(P^4, \rho ^4)$ was trivial.

The other primary invariant, the action of $W$ on $\pi _1(P^4_0, \rho ^4)$, is
trivial as remarked above. To restate the argument, the map $a'$ induces an
injection on $\pi _1$, therefore the $\pi _1$ of the fiber is the image of $\pi
_2$ of the base; but since the base is an infinite loop space, the action of
$\pi _1$ of the base on $\pi _2$ of the base is trivial; thus the action of $\pi
_1$ of the base and in particular of $W$ on  $\pi _1(P^4_0, \rho ^4)$ is
trivial.

We can now look at the secondary invariant for $(P^4, \rho ^4)$.
It is the map
$$
K\rightarrow J^2 = K(H^1_{Dol}(A), 2),
$$
which is the next K\"unneth component of the pullback of $h$ to
$A_{UDol}\times K$,
$$
\bigwedge ^2W \rightarrow H^1(A, \Oo _A)\subset H^1_{Dol} (A) = \pi _2(Hom
(A_{UDol}, K(\Oo , 3))).
$$

This claim tells us that the secondary class for $P^4$ is just $h$ considered
as an element of
$$
H^1(A, \Omega ^2_A)= H^1(A, \bigwedge ^2W^{\ast})= H^1(A,\Oo _A)\otimes _{\Oo}
\bigwedge ^2W^{\ast}.
$$

Recalling that
$$
P^4_0 = \ast \times _{Hom (A_{UDol}, K(\Oo , 3))}K(H^1(A, \Oo _A), 1)
$$
and $\rho ^4= (0,0)$,
we have that $P^4_0$ is just the homotopy fiber of the second morphism.
The long exact sequence for the fibration gives
$$
0\rightarrow \pi _2(Hom (A_{UDol}, K(\Oo , 3)))\rightarrow
\pi _1(P^4_0, \rho ^4) \rightarrow 0
$$
(the morphism from $\pi _1$ of the total space to $\pi _1$ of the base being
injective, and $\pi _2$ of the total space being zero).
This long exact sequence persists under the maps
$$
P^4_0\rightarrow P^3_0\rightarrow P^2_0\rightarrow P^1_0
$$
and furthermore, even into $P_0$ where the long exact sequence for the fiber of
a morphism is replaced by a long exact sequence for the homotopy fiber product.
It follows that the secondary class for $(P, \rho )$ is the image of our above
class (basically $h$)
$$
\bigwedge ^2W \rightarrow H^1(A, \Oo _A)
$$
under composition with the map
$$
H^1(A, \Oo _A)=\pi _2(Hom (A_{UDol}, K(\Oo , 3)), 0)
\rightarrow
\pi _2(Hom (A_{UDol}, T), \rho |_{A_{UDol}})
$$
$$
\rightarrow
\pi _1(Hom (X_{UDol}, T), \rho )=H^2_{Dol}(X)/(\eta ).
$$
One slight twist to notice is that the morphism
$K(\Oo , 3)\rightarrow T$ in question (over $A_{UDol}$) is shifted
by $h$.  This shift is recovered in the last equality,
where we undo a shift by $\varphi$.

This morphism from the long exact sequence for the fiber product (of $Hom$'s)
is equal to the connecting morphism
$$
H^1(A, \Oo _A)\rightarrow H^2(X, \Oo _X) \subset H^2_{Dol}(X)/(\eta ).
$$

Finally, we have concluded that our secondary class is the composition
$$
\bigwedge ^2W^{\ast} \rightarrow H^1(A, \Oo _A) \rightarrow H^2(X, \Oo _X)
$$
where the first map is $h$ (which depends on $\eta$ and which we investigate
below) and the second map is the connecting morphism.

The remaining problem is to calculate
$h\in H^1(A, \Omega ^2_A)$ or more precisely its image by the connecting
morphism.
We have the exact sequence of the cohomology of the pair $(X, A)$ with
coefficients in $\Omega ^2_X$:
$$
0= H^1(X, \Omega ^2_X) \rightarrow
H^1(A, \Omega ^2_A) \rightarrow
$$
$$
H^2((X, A), \Omega ^2_X) \rightarrow
H^2(X, \Omega ^2_X).
$$
The class $\eta$ may be considered as lying in
$H^1((X, A), \Omega ^1_X)$, which is the statement that our map to $T$ factors,
over $A$, through a map to $K(\Oo , 3)$. Therefore the
cup product $\eta \cup \eta$ can be
considered as lying in $H^2((X, A), \Omega ^2_X)$ and mapping to  zero in
$H^2(X, \Omega ^2_X)$.  The class  $h\in H^1(A, \Omega ^2_A)$ is the preimage of
$\eta \cup \eta$ in the exact sequence for the pair $(X,A)$.

Now recall that $X = B\cup C$. This means that the pair $(X,A)$ decomposes
as a ``disjoint union'' of the pairs $(X,B)$ and $(X,C)$ (after applying
excision).  Thus we can write
$$
H^2((X, A), \Omega ^2_X) =
H^2((X, B), \Omega ^2_X) \oplus H^2((X, C), \Omega ^2_X) .
$$

Our class $h\in H^1(A, \Omega ^2_A)$ corresponds
to $b+c$ with
$$
b= -n^2 [H]^2\in H^2((X, B), \Omega ^2_X)
$$
and
$$
c= m^2 [E]^2\in H^2((X, C), \Omega ^2_X).
$$

We now do the same thing for $Z= N\cup C$. The class $h$ corresponds here to an
element of
$$
H^2((Z, N), \Omega ^2_Z) \oplus H^2((Z, C), \Omega ^2_Z) .
$$
This is again of the form $b'+c'$ but now with
$$
b'= -n^2 [H]^2\in H^2((Z, N), \Omega ^2_Z)
$$
and
$$
c'= m^2 [P] \in H^2((Z, C), \Omega ^2_Z).
$$

From the long exact sequence of the pair $(Z, N)$ and the fact that $N$ is
affine we find that
$$
H^2((Z, N), \Omega ^2_Z) \cong H^2(Z, \Omega ^2_Z)=\cc .
$$

To calculate our secondary class we have to contract $h$ with $\alpha \wedge
\beta$ to obtain a class in $H^1(A, \Oo _A)$ and then take its image by the
connecting map in $H^2(X,\Oo _X)$ or equivalently in $H^2(Z, \Oo _Z)$.
To measure this image, use Serre duality: we will  choose a form $\omega \in
H^0(Z, \Omega ^2_Z)$ and take the cup-product to end up with a class in $H^2(Z,
\Omega ^2_Z)$ (of which we then take the residue to end up in $\cc$).
The result of this procedure is the same as if we first contract
$\omega$ with $\alpha \wedge \beta$ and then multiply this section of
$H^0(A, \Oo _A)$ by $h$ getting a class in $H^1(A, \Omega ^2_A)$. Then take the
image of this class by the connecting map and take its residue.
Note that the contraction of $\omega$ with $\alpha \wedge \beta$ is
defined over all of the neighborhood $N$. Call this section $\lambda \in \Oo
(N)$. We are now reduced to the problem of calculating the
image
in $H^2(Z, \Omega ^2_Z)$ under the connecting map for the covering $Z=N\cup C$,
of $\lambda h \in H^1(A, \Omega ^2_A)$.

We have written that the image of $h$ in $H^2((Z, A), \Omega ^2_Z)$
(which we shall denote $[h]$)
is equal to $b'+c'$ where $b'\in H^2((Z, N), \Omega ^2_Z)$ and
$c'\in H^2((Z, C),\Omega ^2_Z)$. The components $b'$ and $c'$ are obtained by
residue maps for the class $h$, along respectively $H$ and $P$
(noting that $H= Z-N$ and $P= Z-C$).  The form of the residue map is not
important for us, just the fact that the classes have poles of order $1$; it
follows that if $\lambda$ is a regular function on $N$ (a neighborhood of $P$)
then the residue of $\lambda h$ at $P$, is equal to $\lambda $ times the
residue of $h$ at $P$. Therefore we can write
$$
[\lambda h] = b'' + \lambda c',
$$
with
$$
b''\in H^2((Z,N), \Omega ^2_Z)
$$
and
$$
\lambda c'\in H^2((Z,C), \Omega ^2_Z)
$$
both being obtained by excision.

The value of $b''\in H^2((Z,N), \Omega ^2_Z)\cong \cc$ is determined by the
condition that the image of $[\lambda h]$ in $H^2(Z, \Omega ^2_Z)$ be
zero (from the long exact sequence for the pair $(Z,A)$).

We now look at the image of $\lambda h$ by the connecting map
in the long exact sequence of the covering $Z= N\cup C$,
$$
H^1(N, \Omega ^2_N)\oplus H^1(C, \Omega ^2_C) \rightarrow
H^1(A, \Omega ^2_A)
$$
$$
\rightarrow H^2(Z, \Omega ^2_Z)\rightarrow \ldots .
$$
We can decompose this connecting map as a composition
$$
H^1(A, \Omega ^2_A) \rightarrow H^2((Z, A), \Omega ^2_Z)
\rightarrow H^2((Z, C), \Omega ^2_Z) \rightarrow H^2(Z, \Omega ^2_Z),
$$
where the second arrow is the projection onto the first factor in the excision
decomposition
$$
H^2((Z, A), \Omega ^2_Z) = H^2((Z, C), \Omega ^2_Z) \oplus
H^2((Z, N), \Omega ^2_Z)
$$
One could equally well use the second factor, with a sign change; our
calculations are not accurate insofar as signs are concerned.

Thus the image of $\lambda h$ by the connecting map
for the covering $Z= N\cup C$, is equal to the class of either $\lambda c'$ or
of $-b''$. We don't know how to calculate $b''$ so we use the representation as
$\lambda c'$. This image is then equal to
$$
\lambda \cdot (m^2 [P])
$$
which is just $m^2\lambda (P)$ (because as noted above, $[P]$ is represented
by cocycles with poles of order $1$).

We have established the formula that the image of
$\alpha \wedge \beta $ under the map
$$
\bigwedge ^2W \rightarrow H^2(Z, \Oo _Z)
$$
is a class which, when paired with a form $\omega \in H^0(Z, \Omega ^2_Z)$,
gives
$$
m^2\omega (\alpha \wedge \beta )(P).
$$
This completes the proof of Theorem \ref{calculation}.
\eop

\numero{APPENDIX I: Relative Breen calculations in
characteristic $0$}

Crucial to the reasonable working of a theory of nonabelian cohomology is the
possibility of calculating the invariants in the Postnikov tower of the spaces
which measure the ``shape''.  In our setup, this means that we would like to
calculate $H^i(K(\Oo , m), \Oo )$.  This calculation is the algebraic analogue
of the classical Eilenberg-MacLane calculations. The algebraic version is the
subject of Breen's work \cite{Breen1}, \cite{Breen2}. His motivation came
mostly from arithmetic geometry, so he concentrated on the case of
characteristic $p$ in \cite{Breen2}. The characteristic $0$ version, while
not explicitly stated in \cite{Breen2}, is implicit there because it is strictly
easier than the characteristic $p$ case: there are no new classes
coming from Frobenius.

These calculations for
the case of base scheme $Spec (\cc )$ are sufficient for the purposes of the
present paper, but eventually a relative version will also be useful. In the
context of calculation of $Ext$ sheaves (i.e. the stable part of the
calculation) this relative version was already evoked in \cite{Breen1}, where
Breen states that the $Ext ^i(G,\cdot )$ for representable group schemes
$G$ can always be calculated.  This part of the topic was not really taken up
afterward, probably for lack of a reasonable category of sheaves over a base
scheme $S$.

Such a category of sheaves is provided by Hirschowitz's notion of {\em
$U$-coherent sheaf} \cite{Hirschowitz}, see also Jaffe's recent paper
\cite{Jaffe}.
\footnote{
Hirschowitz's notion is similar to, but not quite the same as Auslander's
theory of ``coherent functors'' developed in the 1960's. Jaffe's paper
\cite{Jaffe} views $U$-coherent sheaves as a modification or generalization of
Auslander's theory---one looks at functors of algebras rather than functors of
modules. Jaffe, who seems to have been unaware of Hirschowitz's paper, cites
an unpublished letter from Artin to Grothendieck, dating from the 1960's, as a
reference for the generalized version of Auslander's theory. In order to
straighten out the history of this notion, one would have to compare Artin's
letter with \cite{Hirschowitz}.}
We change Hirschowitz's
notation and call these objects {\em vector sheaves}. The category of vector
sheaves over a base $S$ is defined in \cite{Hirschowitz} as the smallest abelian
subcategory of sheaves of $\Oo$-modules on the big site $Sch /S$, containing
$\Oo$ and stable under localization of the base. Thus, locally over $S$ vector
sheaves are obtained starting with $\Oo$ by repeated applications of taking
direct sums, kernels and cokernels.  The abelian category of vector sheaves has
several nice properties \cite{Hirschowitz}. The coherent sheaves on $S$, which
are defined as cokernels  $$ \Oo ^a \rightarrow \Oo ^b \rightarrow \Ff
\rightarrow 0, $$
are vector sheaves. Coherent sheaves are injective objects. Their duals,
which we
call {\em vector schemes}, are the group-schemes with vector space structure
(but not necessarily flat) over $S$. These admit dual presentations as
kernels of maps $\Oo ^b\rightarrow \Oo ^a$. They are projective objects (at
least
if the base $S$ is affine). If $S$ is affine, then any vector sheaf $U$ admits
resolutions $$ 0\rightarrow V\rightarrow V' \rightarrow V'' \rightarrow U
\rightarrow 0 $$
with $V$, $V'$ and $V''$ vector schemes; and
$$
0\rightarrow U \rightarrow \Ff \rightarrow \Ff '\rightarrow \Ff '' \rightarrow 0
$$
with $\Ff$, $\Ff '$, $\Ff ''$ coherent sheaves.

{\bf Example:} The motivating example for the definition of vector sheaf in
\cite{Hirschowitz} was the following example. If $E^{\cdot}$ is a complex of
vector bundles over $S$ then the cohomology sheaves defined on the big site
$Sch /S$ are vector sheaves. Indeed they are vector sheaves of a special type
which Hirschowitz calls
``cohomologies'': quotients of vector schemes by coherent sheaves. This example
is important because it arises from the cohomology of flat families of coherent
sheaves parametrized by $S$: if $f:X\rightarrow S$ is a projective morphism and
$\Ff$ is a coherent sheaf on $X$ flat over $S$ then a classical result says
that the higher direct image complex $R^{\cdot} f_{\ast}(\Ff )$ (calculated on
the big site) is quasiisomorphic to a complex of vector bundles. The notion of
vector sheaf (``$U$-coherent sheaf'' in \cite{Hirschowitz}) thus keeps track of
the jumping of cohomology of flat families of cohoerent sheaves.
\footnote{
This type of example comes up in relation with Dolbeault cohomology: for example
let $$
M:= {\cal M}_{Dol}(X,G) = Hom (X_{Dol}, K(G,1))
$$
be the moduli stack of principal
Higgs $G$-bundles. If $V$ is a representation of $G$ then we obtain
$$
T:= K(V/G, n) \rightarrow K(G,1)
$$
with fiber $K(V,n)$. There is a universal local system $E$ on $X_{Dol}\times
M$, and
$$
\pi _i(Hom (X_{Dol}, T)/M, 0) = H^i(X_{Dol}\times M/M, E).
$$
This is calculated by a Dolbeault complex for $E$ on $X/M$, and the general
discussion of cohomology in flat families applies. Therefore the
$H^i(X_{Dol}\times M/M, E)$ are vector sheaves over $M$ (here $M$ is an
algebraic stack; the condition of being a vector sheaf means that the
pullback to any scheme $Y\rightarrow M$ is a vector sheaf on $Y$). }

The most surprising property from \cite{Hirschowitz} is
that the duality functor $U^{\ast} := Hom (U,\Oo )$ is exact, and is an
involution. This is due to the fact that we take the big site $Sch /S$ rather
than the small Zariski or etale sites.

Another interesting point is that there are two different types of tensor
products of vector sheaves: the {\em tensor product}
$$
U\otimes _{\Oo} V := Hom (U, V^{\ast} )^{\ast},
$$
and the {\em cotensor product}
$$
U\otimes ^{\Oo} V := Hom (U^{\ast}, V).
$$
These are not the same (although they coincide for coherent sheaves cf Lemma
\ref{cohtensor} below) and in particular they don't have the same exactness
properties. Neither of them is equal to the tensor product of sheaves of
$\Oo$-modules. See \cite{kobe} and \cite{RelativeLie} for further discussion.

The above facts work in any characteristic and depend on the
$\Oo$-module structure. However, in characteristic zero vector sheaves have the
additional property that the morphisms $U\rightarrow V$ of sheaves of abelian
groups over $Sch /\cc $ are automatically morphisms of $\Oo$-modules, see
\cite{kobe} \cite{RelativeLie}. Similarly, extensions of sheaves of abelian
groups, between two vector sheaves, are again vector sheaves. These properties
persist for the higher $Ext^i$, see Corollary \ref{ext} below. These properties
do not remain true in characteristic $p$, as shown precisely by Breen's
calculations of \cite{Breen2}. The basic problem is that Frobenius provides a
morphism $\Oo \rightarrow \Oo$ of sheaves of abelian groups, which is not a
morphism of sheaves of $\Oo$-modules. This difficulty in characteristic $p$
seems to be the main obstacle to realizing a reasonable analogue of rational
homotopy theory, for homotopy in characteristic $p$. So we stick to
characteristic $0$!

We don't give a detailed introduction to vector sheaves, rather we refer the
reader to \cite{Hirschowitz}, \cite{kobe} and \cite{RelativeLie}. However, we
do take this opportunity to correct an omission from
\cite{kobe} and \cite{RelativeLie}.
Without the following lemma, the discussion in those references often seems
contradictory, as tensor products and cotensor products are interchanged when
the coefficients are coherent sheaeves.
For example,
in the statement of Corollary 3.9 of \cite{kobe} (which we restate as Theorem
\ref{bc} and prove in more detail below), all terms occuring are coherent
sheaves. Thus the tensor product which appears in the notation is also
equal to the cotensor product. It is the cotensor product which appears
most naturally in that situation. Indeed, Lemma \ref{cohtensor} below explains
(i.e. justifies) the seemingly erroneous statement  $\Ff \otimes _{\Oo} \Gg
=\underline{Hom}(\Ff ^{\ast}, \Gg )$ in the proof of Corollary 3.9 of
\cite{kobe}.

\begin{lemma}
\mylabel{cohtensor}
Suppose $\Ff$ and $\Gg$ are coherent sheaves. Then the tensor product $\Ff
\otimes _{\Oo} \Gg$ and the cotensor product $\Ff \otimes ^{\Oo}\Gg$ coincide.
\end{lemma}
{\em Proof:}
Choose a presentation
$$
\Oo ^a \rightarrow \Oo ^b \rightarrow \Gg \rightarrow 0.
$$
Now $\Ff ^{\ast} := Hom (\Ff , \Oo )$ is a vector scheme, in particular it is a
scheme affine over $S$. Therefore the functor
$$
U\mapsto Hom (\Ff ^{\ast} , U)
$$
is exact in $U$. Applying this functor to the above presentation we obtain
$$
\Ff ^a \rightarrow \Ff ^b \rightarrow Hom (\Ff ^{\ast} , \Gg )\rightarrow 0.
$$
The term $Hom (\Ff ^{\ast} , \Gg )$ is by definition the cotensor product $\Ff
\otimes ^{\Oo} \Gg$. 

Taking the dual of the above presentation we obtain
$$
0\rightarrow \Gg ^{\ast} \rightarrow \Oo ^b \rightarrow \Oo ^a.
$$
Applying the functor
$$
U\mapsto Hom (\Ff , U)
$$
which is exact on the left, we obtain
$$
0\rightarrow Hom (\Ff , \Gg ^{\ast})\rightarrow (\Ff ^{\ast})^a
\rightarrow (\Ff ^{\ast})^b.
$$
Taking the dual we get
$$
\Ff ^a \rightarrow \Ff ^b \rightarrow Hom (\Ff  , \Gg ^{\ast})^{\ast}\rightarrow
0.
$$
This time the term $Hom (\Ff  , \Gg ^{\ast})^{\ast}$ is by definition the
tensor product $\Ff \otimes _{\Oo} \Gg$.

In general there is a natural morphism
$$
U\otimes _{\Oo} V = Hom (U,V^{\ast})^{\ast} \rightarrow Hom (U^{\ast}, V)=
U\otimes ^{\Oo} V.
$$
To define this map we define a trilinear morphism
$$
Hom (U,V^{\ast})^{\ast}\times U^{\ast} \times V^{\ast} \rightarrow \Oo ,
$$
by $(\lambda , \mu , \nu )\mapsto \lambda (\nu \mu )$, the product
$\nu \mu$ being the composed morphism
$$
U\stackrel{\mu}{\rightarrow} \Oo \stackrel{\nu}{\rightarrow} V^{\ast}.
$$
This trilinear map gives a morphism
$$
Hom (U,V^{\ast})^{\ast}\rightarrow Hom (U^{\ast}, (V^{\ast})^{\ast})
$$
then note that $(V^{\ast})^{\ast}=V$.

The above presentations for $\Ff \otimes ^{\Oo} \Gg$ and
$\Ff \otimes _{\Oo} \Gg$ (which are the same) are compatible with this natural
morphism, so the natural morphism is an isomorphism
$$
\Ff \otimes ^{\Oo} \Gg\cong \Ff \otimes ^{\Oo} \Gg .
$$
\eop

We now come to the statement of the ``relative Breen calculations in
characteristic $0$''. The case $S=Spec (k)$ for $k$ a field is due to
\cite{Breen2}, and the relative case for a group scheme was suggested in
\cite{Breen1} (where the case of cohomology with coefficients in the
multiplicative group scheme was treated).

\begin{theorem}
\mylabel{bc}
Suppose $S$ is  a scheme over $Spec (\qq  )$.
Suppose $V$ is a vector scheme over $X$ and suppose $\Ff$ is a coherent sheaf
over $S$. Then for $n$ odd we have
$$
H^i(K(V/S, n)/S, \Ff ) = \Ff \otimes _{\Oo} \bigwedge _{\Oo}^{i/n} (V^{\ast}).
$$
For $n$ even we have
$$
H^i(K(V/S, n)/S, \Ff ) = \Ff \otimes _{\Oo} Sym _{\Oo}^{i/n} (V^{\ast}).
$$
In both cases the answer is $0$ of $i/n$ is not an integer. The multiplicative
structures on the left sides, in the case $\Ff = \Oo$, coincide with the
obvious ones on the right sides. In the case of arbitrary $\Ff$, the
natural structures
of modules over the cohomology with coefficients in $\Oo$, on both sides,
coincide. \end{theorem}

We first recall the following result.

\begin{proposition}
\mylabel{somecomplexes}
Suppose $S$ is a scheme and $V$ is a vector sheaf on $S$. Then the complexes
$$
\ldots \Ff \otimes _{\Oo}\bigwedge _{\Oo} ^j \otimes _{\Oo} Sym _{\Oo}^k V
\rightarrow \Ff \otimes _{\Oo}\bigwedge _{\Oo} ^{j-1}
\otimes _{\Oo} Sym _{\Oo}^{k+1} V
\ldots $$
and
$$
\ldots \Ff \otimes _{\Oo}\bigwedge _{\Oo} ^j
\otimes _{\Oo} Sym _{\Oo}^k V \rightarrow
\Ff \otimes _{\Oo}\bigwedge _{\Oo} ^{j+1}
\otimes _{\Oo} Sym _{\Oo}^{k-1} V \ldots
$$
are exact as sequences of vector sheaves (i.e. as sequences of sheaves on the
site $Sch /\cc $).
\end{proposition}
{\em Proof:} For $\Ff = \Oo$ this is Proposition 3.8 of \cite{kobe}. The proof
is easy, obtained by taking the graded symmetric powers of the cohomologically
trivial complex $V\rightarrow V$ (placing this complex starting in odd or even
degrees, leads to the two cases of the statement).

For a general $\Ff$ note that the cotensor product with a coherent sheaf
is exact---indeed, for any vector sheaf $U$, $Hom (\Ff ^{\ast}, U)$ is exact in
$U$ because $\Ff ^{\ast}$ is represented by a scheme, i.e. an element of the
site $Sch /\cc $. The tensor product is equal to the cotensor product because
both sides are coherent sheaves (here is where we use the hypothesis that $V$
is a vector scheme, i.e. $V^{\ast}$ is a coherent sheaf).
\eop

\subnumero{Proof of Theorem \ref{bc}}

Now we start the proof of Theorem \ref{bc}. Suppose that it is
true for $n\leq m-1$, and we prove it for $n=m$. (The case $m=1$ to start the
induction will be treated separately at the end.) Look at the fiber sequence $$
K(V/S, m-1)\rightarrow \ast \stackrel{p}{\rightarrow} K(V/S, m).
$$
We will look at the Leray spectral sequence for the morphism $p$, for cohomology
with coefficients in $\Ff$.

We may assume that $S$ is affine

Let $A_m$ denote the algebra $H^{\ast}(K(V/S, m)/S, \Oo )$ (this is a sheaf of
algebras over $S$)  and let  $A_m(\Ff )$ denote the $A_m$-module
$H^{\ast}(K(V/S,
m)/S, \Ff )$. Note that if $\Ff$ itself is an algebra-object then $A_m(\Ff )$ is
an $A_m$-algebra (graded-commutative).  By the inductive hypothesis, the $A_k$
and $A_k(\Ff )$ are direct sums of coherent sheaves (coherent in each degree)
for $k\leq m-1$.

The $E_2$ term of our spectral sequence is
$$
H^i(K(V/S, m)/S, H^j(K(V/S, m-1)/S, \Ff ))\Rightarrow H^{i+j}(\ast , \Ff ).
$$
In the case $\Ff = \Oo$ the $E_2$ term has a structure of algebra, and for
arbitrary $\Ff$, a structure of module over that algebra. These are
respectively
$$
A_m(A_{m-1})
$$
and
$$
A_m(A_{m-1}(\Ff )).
$$
By induction we know that
$$
A_{m-1}(\Ff ) = A_{m-1}\otimes _{\Oo} \Ff .
$$

The first possible nonzero differential in the spectral sequence is
$$
H^i(K(V/S, m)/S, H^j(K(V/S, m-1)/S, \Ff ))\rightarrow
$$
$$
H^{i+m}(K(V/S, m)/S, H^{j+1-m}(K(V/S, m-1)/S, \Ff )).
$$
We prove by a second induction on $k$, that for all $\Ff$ the answer is as given
in the theorem, for  $H^i(K(V/S, m)/S,\Ff )$ for all $i\leq k$. Suppose this is
true for $i\leq k-1$. Then the elements of the diagonal complex for the above
differential, ending at $(i,j)= (k, 0)$, are all in the region $i\leq k-1$,
except for the term $(k,0)$. By our second inductive hypothesis (applied to the
cohomology of $K(V/S, m)$ with coefficients in the
coherent sheaves $H^j(K(V/S, m-1), \Ff )$, this complex coincides
with one of the two complexes appearing in Proposition
\ref{somecomplexes}, except maybe for the
last term. However, the complex must be exact at the last stage because
otherwise, what is left over would persist into $E_{\infty}$ contradicting the
answer of the spectral sequence (which must be $\Ff$ in degree $0$ and $0$
otherwise).  Proposition \ref{somecomplexes} gives exactness of the complexes
appearing there. Therefore the $E^{k,0}_2$-term of our spectral sequence must
also coincide with the last term of the complex from Propositon
\ref{somecomplexes}.
For example in the case where $m$ is odd,
the end of the spectral sequence is
$$
\ldots \rightarrow
\Ff \otimes _{\Oo}Sym ^2_{\Oo}(V^{\ast})\otimes_{\Oo}\bigwedge _{\Oo}
^{(k-2m)/m} (V^{\ast})
$$
$$
\rightarrow
\Ff \otimes _{\Oo}V^{\ast} \otimes_{\Oo}\bigwedge _{\Oo} ^{(k-m)/m}
(V^{\ast})
$$
$$
\stackrel{d}{\rightarrow} E^{k,0}_2(\Ff ) \rightarrow 0
$$
(where we denote by $d$ the last differential),
whereas the end of the complex of
Proposition \ref{somecomplexes} is
$$
\ldots \rightarrow
\Ff \otimes _{\Oo}Sym ^2_{\Oo}(V^{\ast})\otimes_{\Oo}\bigwedge _{\Oo}
^{(k-2m)/m} (V^{\ast})
$$
$$
\rightarrow
\Ff \otimes _{\Oo}V^{\ast} \otimes_{\Oo}\bigwedge _{\Oo} ^{(k-m)/m}
(V^{\ast})
$$
$$
\rightarrow  \Ff \otimes _{\Oo}\bigwedge _{\Oo} ^{k/m} (V^{\ast})
\rightarrow 0.
$$
Therefore if $m$ is odd,
$$
E^{k,0}_2 =
\Ff \otimes _{\Oo} \bigwedge _{\Oo}^{k/n} (V^{\ast}).
$$
The same holds  with symmetric power instead of exterior power if $m$ is even.

Cup product gives a bilinear morphism
$$
\mu : E^{k-m,0}_2\times H^m(K(V/S, m), \Oo ) =
E^{k-m,0}_2\times V^{\ast} \rightarrow E^{k,0}_2 .
$$
Let
$$
d':
\Ff \otimes _{\Oo}V^{\ast} \otimes_{\Oo}\bigwedge _{\Oo} ^{(k-2m)/m}
(V^{\ast})
\rightarrow E^{k-m,0}_2(\Ff )
$$
denote the previous differential. We know by induction that $d'$ establishes an
isomorphism between $E^{k-m,0}_2(\Ff ) $ and
$$
\Ff \otimes _{\Oo}\bigwedge _{\Oo} ^{k/m} (V^{\ast})
$$
where this latter is considered as a quotient of the range of $d'$ via
Proposition \ref{somecomplexes}.

Refering to the cup-product morphism $\mu$ considered above and its precursor
$$
\mu ': E^{k-2m,m-1}_2\times V^{\ast} \rightarrow E^{k-m,m-1}_2 ,
$$
we
have the Leibniz formula
$$
d\mu '(a,v) = \mu (d'(a), v)
$$
noting that the term $V^{\ast}$ appearing in the formulas is $E^{m,0}_2$
so the differential acts trivially on the variable $v$.

We have exactly the same formula when the terms $E^{i,j}_2$ are replaced by
their counterparts from the sequences of Proposition \ref{somecomplexes}.
Denote
the multiplication in these counterparts by $\nu$ and $\nu '$ and the
differentials by $\delta$ and $\delta '$.

Call the
isomorphism established by $d$,
$$
\psi ^{k,0}: \Ff \otimes _{\Oo}\bigwedge _{\Oo} ^{k/m} (V^{\ast}) \cong
E^{k,0} _2
$$
and similarly we have isomorphisms (by the inductive hypothesis)
$\psi ^{k-m, 0}$, $\psi ^{k-m, m-1}$ and $\psi ^{k-2m, m-1}$.
The definition of $\psi^{k,0}$ is given by the equation
$$
\psi ^{k,0}(\delta (b))= d \psi ^{k-m, m-1} (b).
$$
Similarly we have
$$
\psi ^{k-m,0}(\delta '(b))= d' \psi ^{k-2m, m-1} (b).
$$
We get
$$
\psi ^{k,0}(\nu (\delta '(a), v))=
\psi ^{k,0}(\delta \nu '(a,v)) = d\psi ^{k-m, m-1} (\nu '(a,v)).
$$
On the other hand, in this region we know by induction that the $\psi$ are
compatible with products. Therefore we get
$$
d\psi ^{k-m, m-1} (\nu '(a,v)) = d \mu '(\psi ^{k-2m,m-1}(a), v)
$$
$$
= \mu (d' \psi ^{k-2m,m-1}(a), v)
$$
$$
=\mu (\psi ^{k-m,0}(\delta '(a)), v).
$$
In all we have established the formula
$$
\psi ^{k,0}(\nu (b, v))
=\mu (\psi ^{k-m,0}(b), v)
$$
for any $b= \delta '(a)$. But $\delta '(a)$ is surjective. This establishes the
compatibility of our isomorphism $\psi ^{k,0}$ with products (given already the
compatibility of $\psi ^{k-m, 0}$).

We note in the above proof that elements of the tensor products are always
(locally on $S$)
finite sums of tensors. This can be seen for example from the proof of Lemma
\ref{cohtensor}. Thus to check compatibility with products, for example, it
suffices to check it on elementary tensors as we have done above.

This completes the proof of the theorem, modulo the case $m=1$ which we now
treat.

We have the fiber sequence
$$
V\rightarrow S \stackrel{p}{\rightarrow} K(V,1).
$$
The higher direct images vanish for coefficients in a coherent sheaf $\Ff$ so
the
Leray spectral sequence implies that
$$
R^ip_{\ast}(\Ff ) = 0, \;\; i> 0
$$
and
$R^0p_{\ast}(\Ff ) $ is a local system on $K(V,1)$ which when restricted to $S$
gives
$$
\Ff [V] \cong \Ff \otimes _{\Oo} Sym ^{\cdot}_{\Oo} (V^{\ast}).
$$
Now we use the complex given in Proposition \ref{somecomplexes}, which is
basically a de Rham complex in our situation:
$$
0\rightarrow \Oo  \rightarrow Sym ^{\cdot}_{\Oo} (V^{\ast})
\ldots \rightarrow Sym ^{\cdot}_{\Oo} (V^{\ast})\otimes _{\Oo}
\bigwedge ^i_{\Oo}
(V^{\ast})\ldots
$$
which can then be tensored by $\Ff$ and remains exact. We can define the
translation action of $V$ on all of the terms, and the exact sequence remains
an exact sequence of sheaves with action of $V$. All of the terms
except for the first one are acyclic by the previous result. Therefore
the cohomology of $K(V,1)$ with coefficients in $\Ff$ is equal to the
cohomology of the complex
$$
\Ff \rightarrow
\ldots \rightarrow \Ff \otimes _{\Oo}
\bigwedge ^i_{\Oo}
(V^{\ast})\ldots .
$$
One can check that the differentials are actually zero, so the cohomology is as
desired.

One should check that the cup-product is equal to the obvious product structure
on the exterior-algebra side of the answer. Instead of doing this (which
as such would seem to
be a difficult task) we proceed as follows. The above construction is functorial
(contravariantly) for morphisms $V\rightarrow V'$. It is easy to see that if one
considers an injection $V\hookrightarrow \Oo ^a$ (which exists locally on $S$
by the definition of vector scheme) then the morphism of functoriality induces a
surjection on cohomology, coming from the surjection $\Oo ^a \rightarrow
V^{\ast}$. Thus to establish a formula for cup-products in cohomology, it
suffices to establish the formula for the case $V=\Oo ^a$. In that case we
can apply the K\"unneth formula, or more precisely remark that the same
K\"unneth formula holds for the cohomology as for the exterior algebra, and that
these formulas are compatible via the above isomorphism. The K\"unneth formulae
are both compatible with cup-products. Thus we can reduce to the case
$V=\Oo$, but here the cohomology is concentrated in degrees $0$ and $1$ so there
are no cup-products to verify (excepting the product with a degree $0$ class
but this case is easy).
\eop

\begin{corollary}
\mylabel{cohVPisVS}
Suppose $S$ is a scheme and $T\rightarrow S$ is a relatively $1$-connected very
presentable $n$-stack. Then for any vector sheaf $V$ on $S$,
$$
H^i(T/S, V)
$$
is a vector sheaf.
\end{corollary}
{\em  Proof:}
We use systematically (without further mention) the fact that the category of
vector sheaves is closed under kernels, cokernels and extensions cf \cite{kobe}
\cite{RelativeLie}. The case of $T= K(U/S,n)$ for $U$ a vector scheme and $V$ a
coherent sheaf is given by Theorem \ref{bc}. The case of coefficients in any
vector sheaf $V$ is obtained by taking a resolution of $V$ by coherent sheaves
and using the long exact sequence of cohomology.  The case of $T=K(U/S,n)$ for
any vector sheaf $U$ is obtained by taking a resolution of $U$ by vector
schemes (divided into two short exact sequences which give rise to two
fibration sequences) and then applying the Leray spectral sequence. Finally, any
relatively $1$-connected very presentable $T$ has a Postnikov tower (relative
to $S$) whose stages are $K(U/S,n)$. Repeated application of the Leray
spectral sequence gives the result.
\eop

\begin{corollary}
\mylabel{homVP}
If $S$ is a scheme and $T\rightarrow S$ and $T'\rightarrow S$ are relatively
$1$-connected very presentable $n$-stacks then $Hom (T/S, T'/S)$ is very
presentable.
\end{corollary}
For this one has to use the fact that $Aut(V)$ is a very presentable group
sheaf when $V$ is a vector sheaf, see \cite{RelativeLie}.
\eop

The following application  was the original motivation for
Breen's calculations of the cohomology of the Eilenberg-MacLane
sheaves \cite{Breen1} \cite{Breen2}. From our version Theorem \ref{bc}, we
obtain
the corresponding result in the relative case in characteristic zero. Similar
corollaries were stated for example for cohomology with coefficients in the
multiplicative group $\Gm$, in \cite{Breen1}.

\begin{corollary}
\mylabel{ext}
{\rm (\cite{kobe} Corollary 3.11)}
Suppose $U,V$ are vector sheaves over a scheme $S$. Let $Ext^i_{\rm gp} (U,V)$
denote the $Ext$ sheaves between $U$ and $V$ considered as sheaves of abelian
groups on $Sch /S$, let $Ext^i_{\rm vs}(U,V)$ denote the $Ext$ sheaves between
$U$ and $V$ considered as vector sheaves on $S$. Then the natural morphisms are
isomorphisms
$$
Ext^i_{\rm vs} (U,V)
\stackrel{\cong}{\rightarrow}Ext^i_{\rm gp} (U,V).
$$
The $Ext^i$ vanish for $i>2$.
\end{corollary}
{\em Proof:}
Let $K_{\cdot}(U,n)$ denote the simplicial presheaf
$$
Y\mapsto K_{\cdot}(U(Y), n )
$$
given by the standard
simplicial Eilenberg-MacLane construction (i.e. Dold-Puppe applied to the
complex with $U$ placed in degree $n$).  We don't take the associated stack
(as doing this or not doesn't affect the morphisms to an $m$-stack). Let
$\zz K_{\cdot}(U,n)$ denote the associated presheaf of simplicial free abelian
groups.  Finally let $N\zz K_{\cdot}(U,n)$ be the presheaf of normalized
complexes (in the homology direction i.e. with differential of degree $-1$) of
this simplicial abelian group. For each $Y$,
$$
N\zz K_{\cdot}(U,n)(Y)
$$
is a complex with homology group $U(Y)$ in degree $n$, and with all other
cohomology groups equal to $0$ in degrees $<2n$. Thus if $\Ff$ is an injective
sheaf of groups then the (cohomological) complex of sheaves
$$
Hom (N\zz K_{\cdot}(U,n), \Ff )
$$
has homology sheaf $Hom (U, \Ff )$ in degree $n$ and zero homology in all other
degrees $< 2n$.
It follows that if $\Ff$ is any sheaf of groups then for $i<n$,
$$
H^{n+i}Hom (N\zz K_{\cdot}(U,n), \Ff ) = Ext ^i_{\rm gp}(U, \Ff ).
$$
On the other hand, this complex of sheaves also calculates
$H^{\cdot}(K(U, n), \Ff )$. Thus we find that
$$
H^{n+i}(K(U,n), \Ff ) = Ext ^i_{\rm gp}(U,\Ff ), \;\; i<n.
$$
This holds true for any sheaves of groups $U$ and $\Ff$. This is one of the
motivating observations of Breen's paper \cite{Breen2}---we have repeated the
proof here for the reader's convenience.

Now suppose that $U$ and $\Ff$ are vector sheaves. If $U$ is a vector scheme
and $\Ff$ is a coherent sheaf then, via the above observation, the relative
Breen
calculations (Theorem \ref{bc}) show that $Ext ^i_{\rm gp}(U,\Ff )=0$ for $i>0$.

A coherent sheaf $\Ff$ is an injective object in the category of vector
sheaves, and the functor $Hom (\cdot , \Ff )$ is exact (cf the discussion of
vector sheaves in \cite{RelativeLie} for example). Thus if $U$ is any vector
sheaf, we can (locally on $S$) resolve it by vector schemes and apply the
previous paragraph. The functor $Ext ^0(\cdot , \Ff )= Hom (\cdot , \Ff
)$ is exact (recall from \cite{kobe} Lemma 3.2 or \cite{RelativeLie} Lemma
4.5 that morphisms of sheaves of abelian groups are the same as morphisms of
vector sheaves so the $Hom$ is the same in the two categories). Using this
exactness we get that $Ext ^i_{\rm gp}(U,\Ff )=0$ for $i> 0$. Finally, if $V$ is
any vector sheaf then we can resolve it by coherent sheaves, which is an
injective resolution in the category of vector sheaves. This is also an acyclic
resolution for $Ext$ in the category of sheaves of abelian groups, so we obtain
the isomorphism
$$
Ext^i_{\rm vs} (U,V)
\stackrel{\cong}{\rightarrow}Ext^i_{\rm gp} (U,V).
$$
The vanishing of the $Ext^i$ for $i>2$ comes from the fact that any vector
sheaf $V$ has a resolution of length $2$ (i.e. with terms in degrees $0,1,2$)
by coherent sheaves (cf \cite{kobe} \cite{RelativeLie}.
\eop

We can apply \ref{ext} to the example discussed at the start of the appendix.
Suppose $S$ is a scheme and suppose $E^{\cdot}$ is a complex of vector bundles
on $S$. The cohomology sheaves $V^i= {\cal H}^i(E^{\cdot})$ are vector sheaves.
In general, a complex with given cohomology objects is determined by higher
extension classes in $Ext^i$ for all values of $i\geq 2$. However, by virtue of
the above theorem the $Ext ^i(V^j, V^k)$ vanish for $i\geq 3$. Thus the complex
$E^{\cdot}$ is determined completely by the successive extension classes
$$
\delta _{j,j+1} \in Ext ^2(V^{j+1}, V^j).
$$
The same is true for any complex of vector sheaves with $V^j$ as cohomology
objects.

{\bf Problem:} describe the conditions which must be satisfied by the classes
 $\delta _{j,j+1}$ for the complex determined by these classes to be
(quasiisomorphic to) a complex of vector bundles.

\numero{APPENDIX II: Representability of very presentable shape}

The following result was stated without proof in
(\cite{kobe}, the discussion above Theorem 5.7).  Since we refer anew to this
result in our discussion after Theorem \ref{calculation} of the present paper, I
felt it to be an opportune time to give a proof.

\begin{theorem}
\mylabel{representable1}
Suppose $\Ff$ is a connected $n$-stack on $Sch /\cc $ such that
the cohomology sheaves $H^i(\Ff , \Oo )$ are represented by finite dimensional
vector spaces.
Suppose furthermore that
$H^0(\Ff , \Oo ) = \Oo $ and $H^1(\Ff , \Oo )= 0$.
Then the $n+1$-functor
$$
T\mapsto Hom (\Ff , T)
$$
from $1$-connected very presentable $n$-stacks of groupoids $T$ to the same,
is representable by a morphism $\Ff \rightarrow \Sigma$, with $\Sigma$ being a
$1$-connected very presentable $n$-stack.
\end{theorem}
{\em Proof:}
It suffices to have a morphism $\Ff \rightarrow \Sigma$
which induces an isomorphism
$$
H^i(\Sigma ,\Oo )\stackrel{\cong}{\rightarrow}
H^i(\Ff ,\Oo )
$$
for any $i$.
We say that a morphism $\Ff \rightarrow \Sigma _m$ is {\em $m$-arranged}
if the induced morphisms on cohomology with coefficients in $\Oo$ are
isomorphisms for $k\leq i < m$ and injective for $k\leq i = m$. Note that the
morphism $\Ff \rightarrow \ast$ is $1$-arranged because of the hypothesis that
$H^1(\Ff , \Oo )=0$. The strategy of the proof (taken from E. Brown
\cite{EBrown}) will be to suppose that we have constructed $\Ff \rightarrow
\Sigma _m$ which is $m$-arranged. Then we will construct a factorization
$$
\Ff \rightarrow \Sigma _{m+1} \rightarrow \Sigma _m
$$
where the first morphism is $m+1$-arranged. By induction this suffices to prove
the theorem (we can stop as soon as we get to $m>n$).

So start with the situation of $\Ff \rightarrow \Sigma _m$, $m$-arranged,
$m\geq 1$. Let
$$
\Cc := Cone (\Ff \rightarrow \Sigma _m)
$$
so we have a map $\Sigma _m \rightarrow \Cc$ which restricted to $\Ff$ gives a
map homotopic to the basepoint $\ast \rightarrow \Cc$ (this basepoint is
included in the definition of $Cone$---it is the vertex of the cone over $\Ff$).

It is easy to see using the $m$-arrangedness
of our map, that the cohomology of $\Cc$ with  coefficients in $\Oo$ vanishes in
degrees $\leq m$.  Furthermore, the $m+1$-st cohomology fits into a long exact
sequence with those of $\Ff $ and $\Sigma _m$:
$$
0\rightarrow H^m(\Sigma _m,\Oo )\rightarrow H^m(\Ff , \Oo )
$$
$$
\rightarrow H^{m+1}(\Cc ,\Oo ) \rightarrow H^{m+1}(\Sigma _m, \Oo )
$$
$$
\rightarrow H^{m+1}(\Ff  , \Oo ) \rightarrow H^{m+2}(\Cc , \Oo )
$$
$$
\rightarrow
H^{m+2}(\Sigma _m, \Oo )\rightarrow \ldots .
$$
The cohomology of $\Sigma _m$ with coefficients in $\Oo$ is a finite
dimensional vector space, by Theorem \ref{bc} (this case is contained in the
original characteristic $0$ version obtainable from \cite{Breen2}). The
property of being represented by a finite dimensional vector space is
closed under extensions, kernels and cokernels (\cite{kobe} Theorem 3.3
and \cite{RelativeLie} Corollary 4.10 and Theorem 4.11). Therefore the
cohomology
of $\Cc$ with coefficients in $\Oo$ is again a (sheaf represented by a) finite
dimensional vector space.

Now let
$$
W^{\ast} := H^{m+1} (\Cc , \Oo )
$$
define the finite dimensional vector space $W$. We get a morphism
$\Cc \rightarrow K(W , m+1)$ which is universal for morphisms to $K(U, m+1)$
with $U$ a finite dimensional vector space.
In particular it induces an isomorphism
$$
H^{m+1}(K(W, m+1), \Oo )\stackrel{\cong}{\rightarrow} H^{m+1}(\Cc , \Oo ).
$$

We will compare the previous long exact sequence with the long exact
sequence that occurs at the start of the Leray-Serre spectral sequence for the
morphism $p:\Sigma _m\rightarrow K(W, m+1)$.
Set
$$
Fib:=Fib(\Sigma _m \rightarrow K(W , m+1) ).
$$
Note that we have a morphism $\Ff \rightarrow Fib$.

The higher direct images
occuring in the Leray-Serre spectral sequence for $p$ are constant local systems
over the base, because $K(W,m+1)$ is $1$-connected.
In other words,
$$
R^ip_{\ast} \Oo = H^i(Fib ,\Oo ) \times K(W,m+1)
\rightarrow K(W,m+1) .
$$
Do a standard type of spectral sequence argument. First of all, for
$k<m$ we prove by induction on $k$ that for all $i\leq k$, the
$H^i(Fib,\Oo ) $ are finite dimensional vector
spaces. Suppose we know this for $k-1$. Then in view of the vanishing of the
cohomology of $K(W, m+1)$ with vector space coefficients (the sheaves
represented
by vector spaces are $\Oo ^a$) the terms $E^{i,j}_2$ with $1\leq i\leq m$ and
$j<k$ vanish; whereas for $j=0$ we have
$$
E^{i,0}_2 = H^i(Fib, \Oo )
$$
because of the fact that $K(W,m+1)$ is $1$-connected. Therefore the terms
$E^{k,0}_2$  persist to $E_{\infty}$ and we have
$$
H^k(Fib, \Oo ) = H^k(\Sigma _m , \Oo )= H^i(\Ff , \Oo ).
$$
In view of the hypothesis of the theorem, this proves the induction step of the
first part of the argument. Incidentally we get that the morphism
$\Ff \rightarrow Fib$ induces an isomorphism on cohomology with coefficients in
$\Oo$ in degrees $k< m$. Now look at the term $E^{m,0}_2= H^m(Fib, \Oo )$. The
only differential concerning it is
$$
d_{m+1}: H^m(Fib , \Oo )\rightarrow H^{m+1}(K(W,m+1) , \Oo )
$$
(noting that we already have $H^0(Fib , \Oo ) = \Oo $).
From our hypothesis which implies that $H^1(Fib, \Oo )=0$ we get, similarly,
that the only differential concerning the term $E^{m+1,0}_2$ is
$$
d_{m+2}: H^{m+1}(Fib , \Oo )\rightarrow H^{m+2}(K(W,m+1) , \Oo ).
$$
From these and the fact that the spectral sequence abuts to the cohomology of
$\Sigma _m$, we get the long exact sequence
$$
0\rightarrow H^m(\Sigma _m,\Oo )\rightarrow H^m(Fib, \Oo )
$$
$$
\rightarrow H^{m+1}(K(W,m+1) ,\Oo ) \rightarrow H^{m+1}(\Sigma _m, \Oo )
\rightarrow
$$
$$
H^{m+1}(Fib , \Oo ) \rightarrow H^{m+2}(K(W,m+1) , \Oo )\rightarrow
H^{m+2}(\Sigma _m, \Oo ).
$$
Remark that $H^{m+2}(K(W, m+1),\Oo )=0$---this comes from Theorem \ref{bc}
and it is here where we use $m\geq 1$.
In particular, the morphism
$$
H^{m+2}(K(W, m+1), \Oo )\rightarrow H^{m+2}(\Cc ,\Oo )
$$
is injective for the trivial reason that the left side is $0$.
Recall that the same induced morphism in degree $m+1$ was an isomorphism
(by the construction of $W$).
Therefore, comparing with the previous long exact sequence and
using the $5$-lemma, we get that the morphism
$$
\Ff \rightarrow Fib
$$
induces isomorphisms on cohomology in degrees $\leq m$ and an injection in
degree $m+1$. In other words this morphism is $m+1$-arranged. Thus we can set
$$
\Sigma _{m+1}:= Fib
$$
and we have completed our inductive construction to prove the theorem.
\eop

{\em Definition:} If $\Ff$ satisfies the condition of
Theorem \ref{representable1} then we obtain the representing $1$-connected
very presentable $\Sigma (\Ff )$ with universal morphism
$$
\Ff \rightarrow \Sigma (\Ff ).
$$
We define
(for any basepoint $f:Y\rightarrow \Ff$)
$$
\pi ^{\rm vp}_i(\Ff \times Y/Y, f):= \pi _i(\Sigma \times Y/Y, f).
$$
In the latter case we usually just take a basepoint $f\in \Ff (Spec \, \cc )$
and then denote this by $\pi ^{\rm vp}_i(\Ff , f)$.

\begin{theorem}
\mylabel{fibration}
Suppose $\Ff$ and $\Gg$ are $n$-stacks with basepoint $g\in \Gg (Spec
(\cc )$, which satisfy the criterion of Theorem \ref{representable1} so their
shapes are representable. Suppose
$$
f:\Ff \rightarrow \Gg
$$
is a morphism of $n$-stacks with the following property (we denote by $\Hh$ the
fiber over $g$):  the local systems $R^if_{\ast} (\Oo )$ are
isomorphic to $\Oo ^{a_i}$ on $\Gg$, and that the morphisms
$$
R^if_{\ast}(\Oo ) |_g \rightarrow H^i(\Hh , \Oo )
$$
are isomorphisms. Suppose that
$$
H^1(\Ff, \Oo ) =H^1(\Gg ,\Oo ) = H^1(\Hh ,\Oo )=\Oo
$$
and
$$
H^1(\Ff, \Oo ) =H^1(\Gg ,\Oo ) = H^1(\Hh ,\Oo )=0.
$$
Then we have a fiber sequence for the representing objects
$$
\Sigma (\Hh )\rightarrow \Sigma (\Ff )\rightarrow \Sigma (\Gg ).
$$
\end{theorem}
{\em Proof:}
The Leray spectral sequence for $f$ is
$$
H^i(\Gg , R^jf_{\ast}(\Oo ))\Rightarrow H^{i+j}(\Ff , \Oo ).
$$
In view of the hypothesis, this becomes
$$
H^i(\Gg , \Oo )\otimes _{ \Oo } H^j(\Hh , \Oo )\Rightarrow H^{i+j}(\Ff , \Oo ).
$$
On the other hand, we obtain a morphism of representing shapes
$$
\Sigma (\Ff )\rightarrow \Sigma (\Gg )
$$
(by the universal property of $\Sigma (\Ff )$).  Let $Fib$ denote the fiber of
$\Sigma (\Ff )$ (over the  image of the point $g$).
Note that the $R^i\Sigma (f)_{\ast} (\Oo )$ are constant on $\Sigma
(\Gg )$ because $\Sigma (\Gg )$ is $1$-connected. In particular
$$
R^if_{\ast}(\Oo ) |_g \stackrel{\cong}{\rightarrow} H^i(Fib , \Oo ).
$$
We obtain the spectral sequence
$$
H^i(\Sigma (\Gg ) , \Oo )\otimes _{\Oo} H^j(Fib , \Oo )
\Rightarrow H^{i+j}
(\Sigma (\Ff ), \Oo ).
$$
Note that the $\pi _i(Fib )$ are finite dimensional vector spaces (using the
long exact sequence of homotopy groups of a fibration, and \cite{RelativeLie}
Theorem 4.11 applied to the case of trivial base $S=\ast $). Thus
$H^j(Fib , \Oo )$ are finite dimensional vector spaces.

The composition
$$
\Hh \rightarrow \Sigma (\Ff )\rightarrow \Sigma (\Gg )
$$
is homotopic to the constant map at the basepoint, so we get a map
$$
\Hh \rightarrow Fib.
$$
This gives  maps $H^i(Fib, \Oo )\rightarrow H^i(\Hh , \Oo )$. We claim that
these
are isomorphisms, which would imply that $\Hh \rightarrow Fib$ is a map
representing the very presentable shape of $\Hh$, in other words $\Sigma (\Hh
)\cong Fib$, thus giving the desired result.

To prove the claim, note that the maps in question are compatible via the
previous identifications, with the maps
$$
R^i\Sigma (\Ff )_{\ast} \Oo \rightarrow R^if_{\ast} (\Oo ).
$$
These in turn fit into a morphism of Leray spectral sequences.  We show using
a spectral sequence argument, by induction on $k$, that for all $i\leq k$ we
have
$$
H^i(Fib, \Oo )\stackrel{\cong}{\rightarrow} H^i(\Hh , \Oo )
$$
or equivalently
$$
R^i\Sigma (\Ff )_{\ast} \Oo \stackrel{\cong}{\rightarrow} R^if_{\ast} (\Oo ).
$$
Suppose this is true for $k-1$. Then look at the term
$E^{0,k}_2 = H^k(\Hh , \Oo )$. When we look at the $r$th differential
$$
d^{0,k}_r:E^{0,k}_r\rightarrow E^{r,k+1-r}_r
$$
the term $E^{r,k+1-r}_r$ has not yet been touched by any term $E^{i,j}$ with
$j\geq k$, and after this differential, the term $E^{r,k+1-r}$ is no longer
touched by any further differentials. We have a morphism of spectral sequences
(the above remarks apply to both) which induces an isomorphism on the
abuttments. It follows that the morphism between spectral sequences induces an
isomorphism on images of $d^{0,k}_r$. Furthermore, the morphism induces an
isomorphism on the intersection (for all $r$) of the kernels of the $d^{0,k}_r$.
This implies that the morphism induces an isomorphism on $E^{0,k}_2$ and we
obtain the inductive step for $k$. This proves the claim and hence the theorem.
\eop

\begin{corollary}
\mylabel{complexify}
Suppose $Y$ is a simply connected finite CW complex. Let $\Ff $
be the constant $n$-stack associated to the constant prestack with values
$Y$ (or
more precisely, with values the $n$-type $\tau _{\leq n}Y$). Fix a basepoint
$y\in Y$ which also gives a basepoint section of $\Ff$. Let $\Ff \rightarrow
\Sigma$ be the morphism representing the shape of $\Ff$. The morphisms
induced by $Y\rightarrow \Sigma (Spec \cc )$,
$$
\pi _i(Y,y)\rightarrow \pi _i(\Sigma , y)
$$
induce isomorphisms
$$
\pi _i(Y,y)\otimes _{\zz} \Oo \cong \pi _i(\Sigma , y).
$$
\end{corollary}
{\em Proof:}
Using the previous theorem, we can reduce by the Postnikov tower to the case $Y=
K(A, n)$ for a finitely generated abelian group $A$. Then the Breen
calculations imply that the morphism
$$
K(A, n)\rightarrow K(A\otimes _{\zz} \Oo , n)
$$
induces an isomorphism on cohomology with coefficients in $\Oo$. This implies
that
$$
\Sigma (K(A,n))= K(A\otimes _{\zz}\Oo , n),
$$
which gives the statement of the corollary.
\eop

{\bf Definition:} Fix $n$. If $Y$ is a $1$-connected finite CW complex, then we
define the {\em complexification of $Y$} denoted $Y\otimes \cc$ to be the
$n$-stack $\Sigma (\Ff )$ representing the very presentable shape of the
constant
$n$-stack $\Ff$ with values $\tau _{\leq n}Y$. Note that this notion depends
on $n$ because we have chosen not to treat the questions arising if we try to
take $n=\infty$.

{\bf Example:}
If we apply this to $Y=S^2$ then we
obtain $\Sigma = S^2\otimes \cc $ as defined in \S 6 above. This is
easy to see because, using the previous theorem, the homotopy sheaves of
$\Sigma$ are $\Oo$ in degrees $2$ and $3$; then there are only two
possibilities for $\Sigma$ and they are distinguished by the vanishing or
nonvanishing of the Whitehead product. As the Whitehead product is nonzero for
$S^2$ and the isomorphisms of the previous theorem are compatible with the
Whitehead product (exercise), this implies that the Whitehead product for
$\Sigma$ is nontrivial, therefore $\Sigma$ is equal to the $T$ defined in \S 6.

We propose the above results as a way of interpreting what it means to look
at the ``complexified homotopy type of a space $Y$''. We could do the same
thing over the ground field $\qq$, and then we propose that this is what it
means to look at the ``rational homotopy type'' of $Y$. This notion is
preserved by base extension of the ground field.

Of course this should all be related to the usual definitions of Quillen,
Sullivan, Morgan, Hain et.al. which refer (excepting Quillen) to algebras of
differential forms. In those theories, base extension is obtained by tensoring
the algebra of forms with the field extension. It has always been somewhat
unclear what geometric interpretation to put on this base-extension process,
and we propose the above theory as a way of obtaining a reasonable
interpretation. We don't, however, get into details of the
relationship between the above theory and the differential-forms theories.

One  advantage of the present formulation is that it explains what is going on
in the non-simply connected case: the shape of the constant
sheaf $\Ff=Const(Y)$ is no longer representable by a very presentable object
(except in fairly restricted cases such as finite fundamental group). Thus, the
object which carries the ``rational homotopy'' information of $Y$ is the shape
itself, rather than the representing object which may not exist. The shape, i.e.
the functor
$$
T\mapsto Hom (\Ff , T)
$$
exists even when $Y$ is not simply connected.

\subnumero{Proof of Theorem \ref{representable0}}

The statement of Theorem \ref{representable0} from \S 2 is very slightly
different from the statement \ref{representable1} given above.  We indicate
here how to get \ref{representable0}. Suppose that $\Ff$ is an $n$-stack on
$Sch /\cc$ such that for any affine algebraic group $G$,
$$
K(G,1)\stackrel{\cong}{\rightarrow} Hom (\Ff , K(G,1)).
$$
In particular this implies that $H^0(\Ff , \Oo )=\Oo $ and $H^1(\Ff , \Oo )=0$.
With the hypothesis that $H^i(\Ff , \Oo )$ are represented by finite
dimensional vector spaces, we can apply Theorem \ref{representable1} to get a
morphism
$$
\Ff \rightarrow \Sigma (\Ff )
$$
universal for morphisms to $1$-connected very prepresentable $T$. Note that
$\Sigma (\Ff )$ is $1$-connected. We have to show that it is also universal for
morphisms to $0$-connected very presentable $T$; suppose that $T$ is one such.
We may choose a basepoint $t$, and let $G= \pi _1(T,t)$ (which is an affine
algebraic group). We have a fiber sequence
$$
T'\rightarrow T \rightarrow K(G,1).
$$
This gives a diagram
$$
\begin{array}{ccccc}
Hom (\Sigma (\Ff ), T')& \rightarrow &Hom (\Sigma (\Ff ), T)&\rightarrow &Hom
(\Sigma (\Ff ), K(G,1))\\
\downarrow & & \downarrow && \downarrow \\
Hom (\Ff , T')&\rightarrow & Hom (\Ff ,T) & \rightarrow & Hom (\Ff , K(G,1)),
\end{array}
$$
where the horizontal sequences are fiber sequences. Since $\Sigma (\Ff )$ is
$1$-connected we have $$
K(G, 1) \stackrel{\cong}{\rightarrow } Hom (\Sigma (\Ff ), K(G,1)),
$$
and the same holds for $\Ff$ by hypothesis. Therefore the vertical map on the
right is an equivalence between $K(G,1)$ and our diagram becomes
$$
\begin{array}{ccccc}
Hom (\Sigma (\Ff ), T')& \rightarrow &Hom (\Sigma (\Ff ), T)&\rightarrow
&K(G,1) \\
\downarrow & & \downarrow && \downarrow =\\
Hom (\Ff , T')&\rightarrow & Hom (\Ff ,T) & \rightarrow & K(G,1).
\end{array}
$$
Now note that $T'$ is a $1$-connected very presentable $n$-stack, so the
vertical arrow on the left is an equivalence. Since the base $K(G,1)$ is
$0$-connected, we can use the long exact sequences of homotopy for these
fibrations to conclude that the vertical morphism in the middle is an
equivalence. This is what we needed to know to establish the universal property
of $\Ff \rightarrow \Sigma (\Ff )$ for Theorem \ref{representable0}.
\eop

\subnumero{A relative version}

While we are on the subject, we give a relative version of Theorem
\ref{representable1}.
Recall \cite{kobe} \cite{RelativeLie} that if $Y$ is a scheme then a
$1$-connected $n$-stack $\Ff \rightarrow Y$ (which can also be thought of as a
$1$-connected $n$-stack on the site $Sch /Y$ of schemes over $Y$) is said to
be {\em very presentable} if for any basepoint section $f: Y'\rightarrow \Ff$
for a scheme $Y'\rightarrow Y'$,
the homotopy group sheaves $\pi _i(\Ff |_{Y'}, f)$ are vector sheaves on $Y'$.
Since (for the present discussion) we have assumed $\Ff$ to be relatively
$1$-connected, the homotopy group sheaves don't depend on the choice of
basepoint (indeed, the choice of basepoint is locally unique up
to homotopy which itself is unique up to---nonunique---homotopy).
Therefore they descend to sheaves of abelian groups $\pi _i(\Ff /Y)$ on $Y$.
For $\Ff$ to be very presentable, it is equivalent to require that these be
vector sheaves.

We introduce the following terminology. We say that a covariant endofunctor $F$
from the category of vector sheaves on $Y$ to itself, is {\em anchored} if
the natural map
$$
F(U)\rightarrow Hom (Hom (F(\Oo ),\Oo ), U)
$$
is an isomorphism for any coherent sheaf $U$ (recall that the coherent sheaves
are the injective objects in the category of vector sheaves). The above natural
map comes from the trilinear map
$$
F(U)\times Hom (F(\Oo ), \Oo ) \times Hom (U,\Oo ) \rightarrow \Oo
$$
defined by $(f,g,h)\mapsto g( F(h)(f))$.

\begin{lemma}
\mylabel{anchored1}
(A) If
$$
0\rightarrow F' \rightarrow F \rightarrow F'' \rightarrow 0
$$
is a short exact sequence of natural transformations between covariant
endofuncturs on the category of vector sheaves over $Y$, then if any two of the
three endofunctors is anchored, so is the third.
\newline
(B) If $F$ is an anchored endofunctor which is also left exact, then $F$ is
representable $F(V)= Hom (W, V)$ for a vector sheaf $W=Hom (F(\Oo ), \Oo )$.
\end{lemma}
{\em Proof:}
(A) follows from the $5$-lemma. For (B) suppose $F$ is a left exact anchored
endofunctor.
Set $W:= Hom (F(\Oo ), \Oo )$. The natural map $F (U)\rightarrow Hom (W, U)$
is an isomorphism for coherent sheaves $U$. On the other hand, both sides are
left exact in $U$. Suppose
$$
0\rightarrow U \rightarrow U' \rightarrow U''
$$
is an exact sequence with $U'$ and $U''$ being coherent sheaves. Then we obtain
exact sequences
$$
0\rightarrow F(U) \rightarrow F(U' )\rightarrow F(U'')
$$
and
$$
0\rightarrow Hom(W,U) \rightarrow Hom(W,U' )\rightarrow Hom(W,U''),
$$
and our natural map is a morphism between these exact sequences inducing
isomorphisms on the last two terms. Thus $F(U)\rightarrow Hom (W,U)$ is an
isomorphism. This completes the proof in view of the fact that (locally on $S$)
any vector sheaf $U$ fits into such a short exact sequence.
\eop

\begin{lemma}
\mylabel{anchored2}
Suppose $V$  is a vector sheaf. Then the endofunctor on the category of vector
sheaves defined by
$$
U\mapsto H^i(K(V,m), U)
$$
is anchored.
\end{lemma}
{\em Proof:}
This follows immediately from Theorem \ref{bc} if $V$ is a vector scheme.
Now suppose that we have an exact sequence
$$
0\rightarrow V' \rightarrow V'' \rightarrow V \rightarrow 0
$$
where $V''$ is a vector scheme, and where we know the lemma for $V'$.
This gives a fibration sequence
$$
K(V'' , m)\rightarrow K(V,m)\rightarrow K(V', m+1),
$$
and taking the cohomology with coefficients in a coherent sheaf $U$ leads to a
Leray spectral sequence
$$
H^i(K(V', m+1), H^j(K(V'', m),U))\Rightarrow H^{i+j}(K(V,m), U).
$$
The cohomology
of the fiber are again coherent sheaves by Theorem \ref{bc}, so by the lemma for
$V'$ the natural map occuring in the definition of ``anchored'' induces an
isomorphism on the $E_2$ terms of the spectral sequence. Since the property of
being anchored is preserved by kernels, cokernels and extensions, we get that
the cohomology of $K(V,m)$ is anchored.
\eop

\begin{corollary}
\mylabel{anchored3}
Suppose $T$ is a relatively $1$-connected very presentable $n$-stack over
a scheme $S$. Then the endofunctor
$$
U\mapsto H^i(T/S, U)
$$
is anchored.
\end{corollary}
{\em Proof:}
Decompose $T$ into a Postnikov tower where the pieces are of the form $K(V,m)$
for vector sheaves $V$
\eop

\begin{theorem}
\mylabel{representable2}
Suppose $S$ is a scheme and $\Ff \rightarrow S$ is a morphism of $n$-stacks on
$Sch /\cc$. Suppose that  for any vector sheaf $V$ over $S$, the cohomology
$H^i(\Ff /S, V)$
is again a vector sheaf over $S$.  Suppose furthermore that
$H^0(\Ff /S, V) = V$ and $H^1(\Ff /S, V)= 0$ for any vector sheaf $V$.
Finally suppose that the functors $V\mapsto H^i(\Ff /S,  V)$ are anchored.
Then
the functor
$$
T\mapsto Hom (\Ff /S , T/S )
$$
from relatively $1$-connected very presentable $n$-stacks of
groupoids $T\rightarrow S$ to the same, is represented by a morphism
$\Ff \rightarrow \Sigma$ over $S$, with $\Sigma \rightarrow S$
being a relatively $1$-connected and very presentable $n$-stack over $S$.
\end{theorem}
{\em Proof:}
Follow the same outline as for the proof of Theorem \ref{representable1}.
We try to find a relatively $1$-connected very presentable $\Sigma
\rightarrow Y$ with a morphism
$$
\Ff \rightarrow \Sigma
$$
inducing an isomorphism on cohomology with coefficients in any coherent sheaf
$U$ on $S$ (the isomorphism for coefficients in any vector sheaf $U$ then follows
by resolving $U$ by coherent sheaves).

We say that a morphism $\Ff \rightarrow \Sigma _m$ is {\em $m$-arranged}
if the induced morphisms on cohomology with coefficients in any coherent sheaf
$U$ on $S$ are isomorphisms for $k\leq i < m$ and injective for $k\leq i = m$.
Note that the morphism $\Ff \rightarrow S$ is $1$-arranged because of the
hypothesis that $H^1(\Ff /S , U )=0$.
Thus we may take $\Sigma _1 := S$.
The strategy of the proof will be to
suppose for some $m\geq 1$ that we have constructed $\Ff \rightarrow \Sigma _m$
which is $m$-arranged. Then we will construct a factorization
$$
\Ff \rightarrow \Sigma _{m+1} \rightarrow \Sigma _m
$$
where the first morphism is $m+1$-arranged. By induction this suffices to prove
the theorem.

Make the same constructions, using the same notations (which we won't repeat
here) as in the proof of Theorem \ref{representable1}. Along the way,
replace the
cohomology with coefficients in $\Oo$ (and the higher direct images of $\Oo$
etc.) by cohomology with coefficients in any coherent sheaf $U$ on $S$.

We obtain the first long exact sequence (actually valid for any vector sheaf
$U$ as coefficients)
$$
0\rightarrow H^m(\Sigma _m/S,U )\rightarrow H^m(\Ff /S, U )
$$
$$
\rightarrow H^{m+1}(\Cc /S,U) \rightarrow H^{m+1}(\Sigma _m/S, U)
$$
$$
\rightarrow H^{m+1}(\Ff  /S, U ) \rightarrow H^{m+2}(\Cc /S, U )
$$
$$
\rightarrow
H^{m+2}(\Sigma _m/S, U )\rightarrow \ldots .
$$
The cohomology of $\Sigma _m$ with coefficients in a vector sheaf is again a
vector sheaf, by Corollary \ref{cohVPisVS} above. The property of being
represented by a finite dimensional vector space is closed under extensions,
kernels and cokernels (\cite{kobe} Theorem 3.3 and \cite{RelativeLie} Corollary
4.10 and Theorem 4.11). Therefore the cohomology of $\Cc$ with coefficients in a
vector sheaf $U$ is again a vector sheaf.

When we come to the construction of $W$ we need to say something
more---this is the reason for introducing the notion of ``anchored'' above.
The functor
$$
U\mapsto H^{m+1}(\Cc , U)
$$
is anchored. This comes from the facts that the cohomology of $\Sigma _m$ is
anchored by Corollary \ref{anchored3}, that the cohomology of $\Ff$ is anchored
by hypothesis, and the fact that being anchored is preserved by kernels,
cokernels and extensions (Lemma \ref{anchored1}).
On the other hand, the fact that the cohomology of $\Cc$ vanishes in degrees
$0<i\leq m$ (note that $m\geq 1$) implies that the above functor is left-exact.
Therefore by Lemma \ref{anchored1} (B) it is representable by a vector
sheaf $W$:
we have
$$
H^{m+1}(\Cc , U) = Hom (W, U).
$$
In particular there is a tautological class in $H^{m+1}(\Cc , W)$
corresponding to a morphism $\Cc \rightarrow K(W, m+1)$, and this morphism is
universal for morphisms from $\Cc$ to things of the form $K(U,m+1)$.
In particular it induces an isomorphism
$$
H^{m+1}(K(W, m+1), U )\stackrel{\cong}{\rightarrow} H^{m+1}(\Cc , U ).
$$
Again set
$$
Fib:=Fib(\Sigma _m \rightarrow K(W , m+1) ).
$$
Note that we have a morphism $\Ff \rightarrow Fib$.

Compare the first long
exact sequence with the long exact sequence that occurs at the start of the
Leray-Serre spectral sequence for the morphism $p:\Sigma _m\rightarrow K(W,
m+1)$,  using the same argument as in the proof of Theorem
\ref{representable1}.  We need to know
that the morphism
$$
H^{m+2}(K(W, m+1), U )\rightarrow H^{m+2}(\Cc ,U )
$$
is injective for cohomology with coefficients in a coherent sheaf $U$ (recall
that only coherent sheaves occur as  coefficients for the cohomology in the
definition of arrangedness---one goes back to the general case after the
induction on $m$ is finished). To prove this we again show that
$H^{m+2}(K(W, m+1), U )=0$ (note that this wouldn't be true if $U$ were not a
coherent sheaf and that is the reason why we restrict to coherent sheaves in
the definition of arrangedness).
In fact, using that $m\geq 1$ and following the argument of
\ref{representable1} we get
that
$$
H^{m+2}(K(W, m+1), U )=Ext ^1(W, U).
$$
However, a coherent sheaf $U$ is an injective object in the category of vector
sheaves \cite{RelativeLie} Lemma 4.17, so $Ext ^1(W, U)=0$. This gives a
proof of
the desired statement. Alternatively one can obtain a proof using a spectral
sequence argument with a resolution
$$
0\rightarrow V \rightarrow V' \rightarrow V'' \rightarrow W \rightarrow 0
$$
of $W$ by vector schemes (decompose this into two short exact sequences and use
a Leray spectral sequence argument for each of the corresponding fibration
sequences).

After that the rest of the argument works exactly the same
way as in Theorem \ref{representable1} (calling upon Theorem \ref{bc} in the
relative case as necessary). We don't repeat this.
\eop

\end{document}